\documentclass[11pt]{article}
\usepackage{graphicx}
\usepackage{amssymb}
\usepackage{amsmath}

\newtheorem{remark}{Remark}[section]
\def\begrem{\begin{remark}\rm}
\def\endrem{\null\hfill\blackbox\end{remark}}

\def\be{\begin{equation}}
\def\ee{\end{equation}}

\author{Thomas Respaud\thanks{IRMA, Université de Strasbourg and INRIA-Nancy-Grand Est, CALVI Project-Team}
\and 
Eric Sonnendr\"{u}cker
  \thanks{IRMA, Universit\'e de Strasbourg and INRIA-Nancy-Grand Est, CALVI Project-Team}
}

\title{Analysis of a new class of Forward Semi-Lagrangian schemes for the 1D Vlasov Poisson Equations} 

\newtheorem{theorem}{Theorem}
\newtheorem{lemma}{Lemma}

\begin{document}
\maketitle

\begin{abstract}
The Vlasov equation is a  kinetic model describing the evolution of a plasma which is a globally neutral gas of charged particles. It is self-consistently  
coupled with Poisson's equation, which rules the evolution of the electric field. 
In this paper, we introduce a new class of  forward Semi-Lagrangian schemes for the Vlasov-Poisson system based on a Cauchy Kovalevsky (CK) procedure for the numerical solution of the characteristic curves.
Exact conservation properties of the first moments of the distribution function for the schemes are derived and
a convergence study is performed that applies as well for the CK scheme as for a more classical Verlet scheme.
A $L^1$ convergence of the schemes will be proved. Error estimates (in  $O(\Delta t ^2+h^2 + \frac{h^2}{\Delta t})$ for Verlet) are obtained, where $\Delta t$ and $h=\max(\Delta x, \Delta v)$ are the discretisation parameters. \end{abstract}


\tableofcontents.


\section{Introduction}
The Vlasov equation describes the dynamics of charged particles in a plasma or in a 
propagating beam. The unknown $f(t, x, v)$ 
which depends on the time $t$, the space $x$ and the velocity $v$ 
represents the distribution function of the studied particles. The 
coupling with the self-consistent electric fields is taken into 
account through the Poisson equation. 

The numerical solution of such systems is most of the time
performed using Particle In Cell (PIC) methods, in which 
the plasma is approximated by macro-particles (see \cite{birdsall}). They are 
advanced in time with the electromagnetic fields which are 
computed on a grid. However, despite their capability to treat 
complex problems, PIC methods are inherently noisy, which becomes 
problematic when low density or highly turbulent regions are studied. 
Hence, numerical methods which discretize the Vlasov equation on a 
grid of the phase space can offer a good alternative to 
PIC methods (see \cite{cheng, filbet1, filbet2, sonnen, vecil}).   
The so-called Eulerian methods can deal with strongly nonlinear 
processes without additional complexity, and are well suited for parallel 
computation (see \cite{virginie}). Moreover, semi-Lagrangian methods 
which have first been introduced in meteorology 
(see \cite{staniforth, ZerWoS05, ZerWoS06}), 
try to take advantage of both 
Lagrangian and Eulerian approaches. Indeed, they allow 
a relatively accurate description of the phase space using a fixed mesh 
and avoid traditional step size restriction using the invariance 
of the distribution function along the trajectories. 

Traditional semi-Lagrangian schemes follow the characteristics backward in time.
In \cite{respaud}, following the idea of Reich  \cite{Reich}, we introduced
a forward Semi-Lagrangian scheme for the Vlasov-Poisson system based on a forward numerical solution of the characteristics using a classical Verlet or Runge-Kutta (order 2 and 4) scheme. 
The Verlet scheme can only be applied for specific differential equations, as for example the characteristics of the Vlasov-Poisson system, but not for more general cases, as the characteristics of the guiding centre or the gyrokinetic approximation of the Vlasov equation. Therefore an alternative to Verlet is necessary. On the other hand, Runge-Kutta schemes, which can be used in the general case, are very costly in our context, especially when going to higher order, as they require a deposition of the charge and the solution of the Poisson equation at intermediate time steps.
We propose here, a new scheme for the characteristics based on a Cauchy-Kowalevsky (CK) procedure, that can be performed up to an arbitrary order. Second and third order are developed in the present paper. We shall also discuss the conservation of the first moment for both Verlet and CK algorithms.

A proof of the convergence of PIC method for the Vlasov-Poisson system was performed by Cottet and Raviart \cite{Cottet}. Proofs of 
convergence and stability of the classical Semi-Lagrangian method applied to the same model were obtained by Besse and Mehrenberger \cite{mehrenberger}. These estimates are made in $L^2$ norm, since $L^{\infty}$ seems out of reach as they explain. They manage to do it because they deal with split methods, and thus only consider constant coefficient transport at each split step. In order to prove convergence in more general cases, the $L^1$ norm seems appropriate, as it enables to use the partition of unity property of the splines. Moreover, 
Despr\'es \cite{Despres1} explains possible advantages of studying $L^1$ convergence instead of more common $L^2$. 

We propose here a proof of $L^1$ convergence of the forward semi-Lagrangian scheme with both
Verlet and CK solution of the characteristics in the particular case of linear spline interpolation.
We also obtain second order error estimates in time and space. 

This paper is organized as follows. In the first part, the continuous problem is presented. In the 
second part, the discrete problem and the numerical scheme to solve it are explained. 
We also prove the exact conservation of the first moment with respect to $v$ at the discrete level for both CK and Verlet schemes. 
Then the convergence of our numerical schemes is proved and finally the schemes are validated and compared on a couple of classical test problems.

\section{The continuous problem}

\subsection{The Vlasov-Poisson model}

Let us consider $f(t,x,v) \geq 0$  the distribution function of positively charged particles in phase-space, and $E(t,x)$ the self consistent electric field. The dimensionless Vlasov Poisson system reads
\begin{equation} \label{eq2.1}
\frac{\partial f}{\partial t}+v\partial_x f + E(t, x)\partial_v f= 0,
\end{equation} 
\begin{equation} \label{eq2.2}
\partial_x E(t, x)= \rho(t,x)=\int_{\mathbb{R}}f(t, x, v)dv - 1, 
\end{equation} 
where $x$ and $v$ are the phase space independent variables. A periodic plasma of period L is considered. So $x \in [0,L]$, $v \in \mathbb{R}$, $t \geq 0$. The functions $f$ and $E$ are submitted to the following conditions
\begin{equation} \label{eq2.3}
f(t,0,v)=f(t,L,v), \forall v \in \mathbb{R} , t\geq 0,
\end{equation}
\begin{equation} \label{eq2.4}
E(t,0)=E(t,L) \Leftrightarrow \frac{1}{L}\int_0^L\int_{\mathbb{R}}f(t,x,v)dvdx = 1, \forall t \geq 0,
\end{equation} 
which translates the global neutrality of the plasma. In order to get a well-posed problem, a zero-mean electrostatic condition has to be added, which corresponds to a periodic electric potential:
\begin{equation} \label{eq2.5}
 \int_0^L E(t,x)dx =0, \quad \forall t \geq 0,
\end{equation} 
 and an initial condition
  \begin{equation} \label{eq2.6}
 f(0,x,v)=f_0(x,v), \quad \forall x \in [0,L], v\in \mathbb{R}.
\end{equation} 
 Assuming that the electric field is smooth enough, equations \eqref{eq2.1}, \eqref{eq2.3} and \eqref{eq2.6} can be solved in the classical sense as follows. 
 
 The first order differential system
  \begin{eqnarray} \label{eq2.7}
 \frac{dX}{dt}(t;(x,v),s)&=&V(t;(x,v),s), \nonumber\\
 \frac{dV}{dt}(t;(x,v),s)&=&E(t,X(t;(x,v),s)), 
 \end{eqnarray} 
 where (X(t;(x,v),s),V(t;(x,v),s)) are the characteristic curves, solutions of \eqref{eq2.7} at time $t$ with the initial condition
  \begin{equation} \label{eq2.8}
 X(s;(x,v),s)=x,
 V(s;(x,v),s)=v.
\end{equation} 
For the existence, the uniqueness and the regularity of the solutions of this differential system, the reader is referred to \cite{bouchut}.
 The solution of problem \eqref{eq2.1}, \eqref{eq2.6} is then given by
  \begin{equation} \label{eq2.9}
 f(t,x,v)=f_0( X(0;(x,v),t),V(0;(x,v),t)), \quad  \forall x \in [0,L], v\in \mathbb{R}, t\geq 0.
 \end{equation} 
 Since 
  $$
 \frac{\partial(X,V)}{\partial(x,v)} =1,
 $$
 the conservation of particles is ensured for all times:
  $$
 \frac{1}{L}\int_0^L\int_\mathbb{R}f(t,x,v)\,dv\,dx= \frac{1}{L}\int_0^L\int_\mathbb{R}f_0(x,v)\,dv\,dx=1.
 $$
 
 According to previous considerations, an equivalent form of the Vlasov-Poisson periodic problem is to find $(f,E)$, smooth enough, periodic with respect to $x$, with period $L$, and solving the equations (2.2), \eqref{eq2.7}, \eqref{eq2.8} and \eqref{eq2.9}.
 Introducing the electrostatic potential $\varphi \equiv \varphi(t,x)$ such that $E(t,x)=-\partial_x\varphi(t,x)$, and setting $G=G(x,y)$ the fundamental solution of the Laplacian operator in one dimension. That is
 $-\partial_x^2G(x,y)=\delta_0(x-y)$ with periodic boundary conditions. It comes
  $$
 E(t,x)=\int_0^LK(x,y)(\int_{\mathbb{R}}f(t,y,v)dv-1)\,dy,
 $$
  where
  \begin{eqnarray*}
 K(x,y)=-\partial_xG(x,y) & = & (\frac{y}{L}-1), ~~~ 0 \leq x < y, \nonumber\\
                                           & = & \frac{y}{L},  ~~~ y < x \leq L.
  \end{eqnarray*}        
  
  \subsection{Existence, uniqueness and regularity of the solution of the continuous problem}
  
\begin{theorem}                               
  Assuming that $f_0 \in W_{c,per_x}^{1,\infty}(\mathbb{R}_x \times \mathbb{R}_v)$ ($W_{c,per_x}^{1,\infty}(\mathbb{R}_x \times \mathbb{R}_v$ being the Sobolev space of functions with first derivatives in $L^{\infty}$, compactly supported in v and periodic in x), positive, periodic with respect to the variable $x$ with period $L$, and $Q(0) \leq R$, with $R>0$ defined as follows
  $$
 Q(t):=1+sup{\vert v \vert: x \in [0,L], \tau \in [0,t] \vert f(\tau,x,v)\neq 0},
 $$
 and
$$
 \frac{1}{L}\int_0^L\int_\mathbb{R}f_0(x,v)\,dv\,dx=1,
$$
then the periodic Vlasov-Poisson system has a unique classical solution $(f,E)$, periodic in $x$, with period $L$, for all t in $[0,T]$, such that
$$
f \in W^{1,\infty}(0,T;W_{c,per_x}^{1,\infty}(\mathbb{R}_x \times \mathbb{R}_v)),
$$
$$
E \in  W^{1,\infty}(0,T;W_{per_x}^{1,\infty}(\mathbb{R})),
$$
and there exists a constant $C=C(R,f_0)$ dependent of $R$ and $f_0$ such that
$$
Q(T)\leq CT.
$$ 
 Moreover if we assume that $f_0 \in  W_{c,per_x}^{m,\infty}(\mathbb{R}_x \times \mathbb{R}_v)$, then $(f,E) \in W^{m,\infty}(0,T; W_{c,per_x}^{m,\infty}(\mathbb{R}_x \times \mathbb{R}_v))\times  W^{1,\infty}(0,T;W_{per_x}^{1,\infty}(\mathbb{R}))$, for all finite time $T$.
\end{theorem}

For the proof, the reader is referred to the references \cite{bouchut, Glassey}.

\section{The discrete problem}

\subsection{Definitions and notations}

Let $\Omega=[0,L[\times[-R,R]$, with $R>Q(T)$, and $M_h$ a cartesian mesh of the phase-space $\Omega$. $M_h$ is given by a first increasing sequence $(x_i)_{i \in [0..N_x]}$ of the interval $[0,L]$ and  a second one $(v_j)_{j \in [0..N_v]}$ of the interval [-R,R]. Let $\Delta x_i = x_{i+1}-x_i$ the physical space cell width and  $\Delta v_j = v_{j+1}-v_j$ the velocity space cell width. In order to simplify the study, a regular mesh will be used, i.e $\Delta x_i=\Delta x=\frac{L}{N_x+1}$, and $\Delta v_j=\Delta v=\frac{2R}{N_v}$, where $N_x$, $N_v$ belong to $\mathbb{N}$. Then $h$ is defined being $\max(\Delta x,\Delta v)$.

For each function $g$ defined on all the points $(x_i,v_j)\in M_h$ we will set $g_{i,j}:=g(x_i,v_j)$, and the sequence is completed on $\mathbb{Z} \times \mathbb{Z}$ by periodicity in $x$ and by 0 in $v$. The sequence $(x_i,v_j)$ will also be defined on the whole set $\mathbb{Z} \times \mathbb{Z}$, by $x_i:= i\Delta x$, and $v_j:=-R+j\Delta v$. The set of all L-periodic functions  in $x$ and compactly supported in $v$ will be denoted $P(\Omega)$.

Now let $f=(f_{i,j})_{(i,j)\in \mathbb{Z} \times \mathbb{Z}}$ be a continuous grid-function, periodic in the $x$ direction and compactly supported in the $v$ direction, with a support included in $[-R,R]$. 
If $f$ is a function defined on the points $M_h$, a discrete grid-function $\tilde{f}$ can be defined by $\tilde{f}_{i,j}:=f(x_i,v_j)$ for all $(i,j) \in [0,N_x] \times [0,N_v]$. In order to lighten notations, $f$ will be kept instead of $\tilde{f}$. Let $L^2_h(\Omega)$ (resp. $L^1_h(\Omega)$), the set of grid-functions whose $\vert \vert . \vert \vert_{L^2_h(\Omega)}$ (resp $\vert \vert . \vert \vert_{L^1_h(\Omega)}$) is bounded
$$
\vert \vert f \vert \vert_{L^2_h(\Omega)} = (\Delta x \Delta v \sum_{i=0}^{N_x}\sum_{j=0}^{N_v} \vert f_{i,j}\vert^2)^{\frac{1}{2}},
$$
$$
\vert \vert f \vert \vert_{L^1_h(\Omega)} = \Delta x \Delta v \sum_{i=0}^{N_x}\sum_{j=0}^{N_v} \vert f_{i,j}\vert.
$$

As was precised in the introduction, a $L^2$ convergence analysis for backward Semi-Lagrangian scheme, in the case of a Strang split time advance, was performed in \cite{mehrenberger}. In paper \cite{respaud}, it is explained that for split methods, where the split steps consist of constant coefficient transport, forward and backward methods are exactly the same. So all the $L^2$ results exposed in \cite{mehrenberger} are  also valid for our method when time splitting is used. In this paper, we shall consider the convergence of a non split method, and the $L^1$ norm seems more appropriate for this kind of study. 

\begin{remark}
If $f \in L^2_h(\Omega)$, then $f \in L^1_h(\Omega)$ thanks to the Cauchy Schwarz inequality, as $\Omega$ is bounded.
\end{remark}

In the sequel, a final time T is fixed, as well as a uniform time discretization $(t^n)_{n \leq N_T}$ of the interval $[0,T]$, with  time step $\Delta t = t^{n+1}-t^n$. At each point $(x_i,v_j) \in M_h$, an approximation $f_h(t^n,x_i,v_j)$ of the exact distribution function $f(t^n,x_i,v_j)$ at time $t^n = n\Delta t$ is defined.
The approximation function $f_h(t^n)$ is then given at each point of $\mathbb{R}_x \times \mathbb{R}_v$ thanks to an interpolation operator $R_h$ defined on a uniform grid:
\begin{align*}
R_h: \; L^1(\Omega) \cap P(\Omega)  &\longrightarrow  L^1(\Omega) \cap P(\Omega),\\
 f  &\mapsto  R_h f = \sum_{(i,j) \in \mathbb{Z} \times \mathbb{Z}} f_{i,j} \Psi_{i,j},
\end{align*}
where $\Psi_{i,j}$ will be linear spline functions for our study. In numerical results, since linear interpolation is quite diffusive, cubic splines will be used. In order to get a convergent scheme, the operator $R_h$ must satisfy some approximation properties which will be detailed later.

\subsection{The numerical scheme}

The electric field operator for the real-valued function $g \in L^1([0,L]\times \mathbb{R})$ is defined this way:
\begin{eqnarray}\label{eq3.10}
E[g](x) = \int_0^L K(x,y)(\int_{\mathbb{R}}g(y,v)dv-1).
\end{eqnarray}
The approximate function $f_h$ is solution on the grid of the following Vlasov equation:
$$
\frac{\partial f_h}{\partial t}(t,x,v)+v\frac{\partial f_h}{\partial x}(t,x,v)+E_h(t,x)\frac{\partial f_h}{\partial v}(t,x,v)=0.
$$
This function follows approximate characteristics, solutions of

\begin{eqnarray}\label{eq3.11}
\frac{dX_h}{dt}(t;(x,v),s)=V_h(t;(x,v),s), \nonumber\\
\frac{dV_h}{dt}(t;(x,v),s)=E_h(t,X(t;(x,v),s)), 
\end{eqnarray} 
where $E_h$ is defined exactly from $f_h$ using \eqref{eq3.10}: $E_h=E[f_h](x)$.
So we get:
$$
\forall t \in [t^n,t^{n+1}[: f_h(t,x,v)=\sum_{k,l} \omega_{k,l}^n S_h(x-X_h(t;(x_k,v_l),t^n))S_h(v-V_h(t,(x_k,v_l),t^n)),
$$
so that $f_h$ is given on the mesh at time $t^n$ by:
$$
f_h(t^n,x_i,v_j)=\sum_{k,l} \omega_{k,l}^n S_h(x_i-x_k)S_h(v_j-v_l) \quad \forall n.
$$
These are the interpolation conditions enabling to define $f_h$ everywhere.
The computation of  $(\omega_{k,l}^n)_{k,l}$ from the grid values amounts to solving a linear system, which is trivial in the case of linear splines, where $\omega_{k,l}^{n}=f_h(t^n,x_k,v_l)$.

Let us recall that
the linear B-spline $S$ is defined as follows 
$$
S(x)=
\left\{
\begin{array}{ll}
(1-|x|) & \textrm{if} \ 0\le |x|\le 1,\\
0& \textrm{otherwise.}
\end{array}
\right.
$$
The spline $S_h$ actually used, which shall be defined with the size of the mesh will be
$$
S_h(x)=S(\frac{x}{h})
$$
with $h=\Delta x$ for splines in the $x$ variable and $h=\Delta v$ for splines in the $v$ variable.
From now on, $S_h$ will be denoted by $S$ for the sake of simplicity.

The distribution function is updated this way: The ending point of the characteristic starting from $(x_i,v_j)$ is computed: $(X_h(t^{n+1};(x_i,v_j),t^n),V_h(t^{n+1};(x_i,v_j),t^n)), \forall (i,j)$. Then, since $f_h$ is constant along the approximate characteristics, the value is deposited on the nearest grid points, the number of which depending on the degree of the splines used for the interpolation.  
This amounts to computing $f_h$ at time $t_{n+1}$ at the grid points using the following formula.
\begin{eqnarray*}
f_h(t^{n+1},x_i,v_j)=\sum_{k,l} f_h(t^n,x_k,v_l)S(x_i-X_h(t^{n+1};(x_k,v_l),t^n))S(v_j-V_h(t^{n+1};(x_k,v_l),t^n)) \quad \forall (i,j) ,
\end{eqnarray*}
Note that
 $$(X_h(t^{n+1};(x_i,v_j),t^n),V_h(t^{n+1};(x_i,v_j),t^n)), \forall (i,j),$$
 are computed by a numerical solution of the differential system \eqref{eq3.11}. Since this requires an explicit solution of that system,  any standard  ODE solver such as Verlet, Runge-Kutta or others can be used. Our analysis will be based on the Verlet algorithm, which is second order accurate, and on a Cauchy Kovalevsky procedure, which can be of any order, as an alternative to more costly Runge-Kutta solvers. But we will consider only the second and third order, because higher ones would not increase accuracy in our case, as we will explain.

\subsubsection{Verlet algorithm}

Starting at time $t^n$ from the grid point $(x_k,v_l)$
\begin{itemize}
\item Step $1$: $\forall k,l$, $x_{k,l}^{n+\frac{1}{2}} -x_k \,=\,  \frac{\Delta t}{2} \, v_l $,
\item Step $2$: compute the electric field at time $t^{n+\frac{1}{2}}$
 \begin{itemize}
\item   deposition of the particles $x_{k,l}^{n+\frac{1}{2}}$ on the spatial grid $x_i$ for the density $\rho_h$:
$\rho_h(x_i, t^{n+\frac{1}{2}})=\sum_{k,l} \omega^n_{k,l} S(x_i-x^{n+\frac{1}{2}}_{k,l})$, like in a PIC method.
\item solve the Poisson equation on the grid $x_i$: $E(x_i, t^{n+\frac{1}{2}})$.   
\end{itemize}  
\item Step $3$:  $\forall k,l$, $v_{k,l}^{n+1} -v_l \,=\, \Delta t  \,E(x_{k,l}^{n+\frac{1}{2}}, t^{n+\frac{1}{2}})$,
\item Step $4$: $\forall k,l$, $x_{k,l}^{n+1} -x_{k,l}^{n+\frac{1}{2}} \,=\, \frac{\Delta t}{2}  \, v_{k,l}^{n+1}$.   
\end{itemize}
This is the way the algorithm is implemented. In our convergence study, the slight difference is that an exact solution of Poisson's equation, based on the Green formula \eqref{eq3.10}, is used.

\subsubsection{Cauchy Kovalevsky procedure}

The idea is to get high order approximations of the characteristics using Taylor expansions in time. And then, using the charge conservation equation, and higher velocity moments of the Vlasov equation, to replace time derivatives with terms containing only spatial derivatives and moments at time $t_{n}$ which can be easily computed.
Up to third order, these Taylor expansions in time lead to
\begin{eqnarray*}
X^{n+1}&=&X^n+\Delta t V^n + \frac{\Delta t^2}{2}E^n(X^n)+ \frac{\Delta t^3}{6}\frac{d}{dt}E(X(t),t)_{\vert t=t^n}.\\
V^{n+1}&=&V^n+\Delta t E^n(X^n) + \frac{\Delta t^2}{2}\frac{d}{dt}E(X(t),t)_{\vert t=t^n} + \frac{\Delta t^3}{6}\frac{d^2}{dt^2}E(X(t),t)_{\vert t=t^n}.
\end{eqnarray*}
In order to be able to compute all terms of these expansions we need the three first total time derivatives of $E(X(t),t)$.
\begin{eqnarray*}
\frac{d}{dt}E(X(t),t) & = & \frac{\partial E}{\partial t}(X(t),t)+\frac{dX}{dt}(t)\frac{\partial E}{\partial x}(X(t),t) \nonumber\\
	& = & - J(X(t),t)+\bar{J}(t)+V(t)\rho(X(t),t),
\end{eqnarray*}
where $\rho(x,t)=\int f(x,v,t)\,dv-1$, $J(x,t)=\int f(x,v,t)v\,dv$ and $\bar{J}(t)=\frac{1}{L}\int_{0}^L J(x,t)\,dx$.
Indeed, the Poisson's equation yields $\frac{\partial E}{\partial x}=\rho$ and integrating the Vlasov equation with respect to velocity, yields the charge conservation equation $\frac{\partial \rho}{\partial t}+
\frac{\partial J}{\partial x}=0$.
Hence taking the derivative of the Poisson's equation with respect to time and using this equation we get
$$\frac{\partial}{\partial x}(\frac{\partial E}{\partial t}+J)=0.$$
From which we obtain, as $\int_{0}^L E(x,t)\,dx=0$, that
$$\frac{\partial E}{\partial t}=-J+\bar{J}.$$
The second order total derivative in time of $E$ reads
\begin{eqnarray*}
\frac{d^2}{dt^2}E(X(t),t) & = & - \frac{\partial J}{\partial t}(X(t),t)-V(t)\frac{\partial J}{\partial x}(X(t),t)+\frac{d\bar{J}}{dt}(t)\nonumber\\
	& +& E(X(t),t)\rho(X(t),t) + V(t)(\frac{\partial \rho}{\partial t}(X(t),t)+V(t)\frac{\partial \rho}{\partial x}(X(t),t)).
\end{eqnarray*}
  In order to use this expression, we need $  \frac{\partial J}{\partial t}, \frac{\partial J}{\partial x}, \frac{\partial \rho}{\partial t},
  \frac{\partial \rho}{\partial x},\frac{d\bar{J}}{dt}(t).$
      
The Cauchy-Kovalevsky procedure consists in getting rid of time derivatives, replacing them with space derivatives obtained from the equation, in our case, we use the velocity moments of the Vlasov equation.
First for $\rho$, we use the charge conservation equation:
 \begin{eqnarray}\label{eq3.12}
\frac{\partial \rho}{\partial t}(X(t),t) = -  \frac{\partial J}{\partial x}(X(t),t) .
\end{eqnarray}
 In order to get the time derivative of the current $J$, we need to use the Vlasov equation, multiply it with $v$, and integrate it with respect to $v$, so that we get:
  \begin{eqnarray*}
 \frac{\partial J}{\partial t} + \frac{ \partial}{\partial x} I_2 + E \int_{\mathbb{R}}\frac{\partial f }{\partial v}vdv=0,
 \end{eqnarray*}
 where $I_n(x,t)=\int_{\mathbb{R}}f(x,v,t)v^n \,dv$
 so that, using that $f$ is compactly supported and integrating by parts:  
 \begin{eqnarray}\label{eq3.13}
 \frac{\partial J}{\partial t} (X(t),t) = -\frac{\partial I_2}{\partial x}(X(t),t) + E(X(t),t)(1+\rho(X(t),t)).
\end{eqnarray}
Let us prove that $\frac{d\bar{J}}{dt}(t)=0$

\begin{eqnarray*}
\frac{d\bar{J}}{dt}(t) & = & \frac{1}{L}\int_0^L\frac{\partial J}{\partial t}(t,x)dx\nonumber\\
 & = & \frac{1}{L}\int_0^L(-\frac{\partial I_2}{\partial x}(t,x) + E(t,x)(1+\rho(t,x))dx\nonumber\\
 & = & \frac{1}{L}([I_2(t,0)-I_2(t,L)]+\int_0^LE(t,x)dx+\int_0^LE(t,x)\rho(t,x)dx\nonumber\\
 & = & \frac{1}{L}[E^2(t,L)-E^2(t,0)]\nonumber\\
 & = & 0.
 \end{eqnarray*}
 thanks to periodicity, in fact \eqref{eq2.3}, \eqref{eq2.5}. We will see later, that numerically this value is also  zero.
We finally get the following third order Cauchy Kovalevsky (CK3) time algorithm, using \eqref{eq3.12}, \eqref{eq3.13}:
 \begin{eqnarray*}
 X^{n+1} & = & X^n+\Delta t V^n + \frac{\Delta t ^2}{2}E^n(X^n)
 	+  \frac{\Delta t ^3}{6}(V^n\rho^n(X^n)-J^n(X^n)+\bar{J}),\\
V^{n+1} & = & V^n+\Delta t E^n(X^n) + \frac{\Delta t ^2}{2}(V^n\rho^n(X^n)-J^n(X^n)+\bar{J})\nonumber\\
	& + &  \frac{\Delta t^3}{6}(\frac{\partial I_2}{\partial x}(X^n,t^n)-E^n(X^n)-2V^n\frac{\partial J}{\partial x}(X^n,t^n) +(V^n)^2\frac{\partial \rho}{\partial x}(X^n,t^n)).
\end{eqnarray*}
Let us now introduce a notation, which will be useful later: 
$$
X^{n+1}=X^n+\Delta t V^n + \frac{\Delta t ^2}{2}E^n(X^n) +  \frac{\Delta t ^3}{6}\phi^n(X^n,V^n), 
$$
and
$$
V^{n+1}=V^n + \Delta t E^n(X^n) +  \frac{\Delta t ^2}{2} \phi^n(X^n,V^n)  + \frac{\Delta t ^3}{6}\varphi^n(X^n,V^n),  
$$
where $\phi, \varphi$ are naturally defined. 

\begin{remark}
Obviously, in order to get a second order algorithm (CK2), we just keep the terms until  $\Delta t ^2$ included.
\end{remark}

\subsection{Exact conservation of number of particles and momentum}

\subsubsection{B-spline interpolation}

First, let us recall some useful properties of B-splines interpolation. The linear space of B-splines of order $m+1$ writes, denoting by $s^{(m)}$ the $m^{th}$ derivative of s
$$
S_{m+1,\Delta_x} = \lbrace s(x) \in C^{m-1}(\mathbb{R}), \quad s^{(m+1)}(x)=0, \forall x \in (x_i,x_{i+1}), \forall i \in \mathbb{R} \rbrace,
$$
if $m+1$ is even, and 
$$S_{m+1,\Delta_x} = \lbrace s(x) \in C^{m-1}(\mathbb{R}), \quad s^{(m+1)}(x)=0, \forall x \in (x_{i-\frac{1}{2}},x_{i+\frac{1}{2}}), \forall i \in \mathbb{R} \rbrace,
$$
if $m+1$ is odd.

The space of B-spline functions in two dimensions is defined as the tensor product of $1D$ spaces.
Let us precise the interpolation operator: 
$$
R_{h_{i,j}}(f)(x,v)=\omega_{i,j}(f)S(x-x_i)S(v-v_j),
$$
$$
R_hf(x,v)=\sum_{i,j}R_{h_{i,j}}(f)(x,v).
$$
Now come the properties:
\begin{itemize}
\item $\mathbb{S}_{m+1,h}=Span(S_{m+1}(.-x_i)S_{m+1}(.-v_j); \forall (i,j) \in \mathbb{Z})$,
\item $\mathbb{S}_{m+1,h} \subset W^{k,p} \quad 1\leq p \leq \infty \quad 0 \leq k \leq m$,
\item Stability $\vert\vert R_hf\vert \vert_{L^p(\Omega)} \leq C \vert\vert f \vert \vert_{L^p(\Omega)} \quad \forall f \in L^p(\Omega) \cap P(\Omega), \quad 1 \leq p \leq \infty$ \quad (i),
\item Consistency and accuracy. There exists $C > 0 \quad \vert \quad  \vert\vert f-R_hf \vert \vert_{W^{k,p}(\Omega)} \leq C h^{m+1-k} \vert f \vert_{W^{m+1,p}(\Omega)} \quad \forall f \in W^{m+1,p}(\Omega)\cap P(\Omega) \quad 1\leq p \leq \infty \quad 0 \leq k \leq m$ \quad (ii),
\item $\sum_i S_m(.-x_i) =1$ \quad (iii), $\quad \int S_m(u)du=\Delta x$ \quad (iv).
\item $\sum_l v_lS_1(v_l-v)=v$.  \quad (v)
\end{itemize}
For the last item, we will give the proof:
Let us suppose  $v = v_p+\alpha \Delta v, \alpha \in [0,1[$
\begin{eqnarray*}
\sum_l v_lS_1(v_l-v) & = & v_pS_1(\alpha \Delta v) + v_{p+1}S_1(\Delta v-\alpha \Delta v)\nonumber\\
 		& = & (p+\alpha \Delta v)\nonumber\\
		& = & v
\end{eqnarray*}
Let us also precise particle and momentum conservation. The proof for the mass is independent from the spline degree, and the one for the first moment will only be shown for linear splines, even though it has been checked for the first three splines.

\subsubsection{Particle conservation}

The discrete algorithm preserves the total number of particles, as the following computation shows:
\begin{eqnarray*}
m^{n+1} &=& \int f_h(t^{n+1},x,v)\,dx\,dv,\nonumber\\
	&=&\sum_{i,j} \omega_{i,j}^{n+1}\int S(x-x_i)S(v-v_j)\,dx\,dv,\nonumber\\
	&=&\Delta x \Delta v\sum_{i,j} f_h^{n+1}(x_i,v_j),\nonumber\\ 
      &=&\Delta x \Delta v\sum_{i,j} \sum_{k,l} \omega_{k,l}^{n} S(x_i-X(t^{n+1}; (x_k,v_l), t^n)S(v_j-V(t^{n+1};(x_k,v_l),t^n),\nonumber\\  
   	&=& \Delta x\Delta v\sum_{k,l} \omega_{k,l}^n = \Delta x\Delta v\sum_{i,j} f^{n}(x_i,v_j)= m^n,
\end{eqnarray*}
thanks to partition of unity property (iii).

Let us precise the way the interpolation operator acts, in fact:
\begin{eqnarray}\label{eq4.15}
f_h(t^{n+1},x,v) & = & R_h(\sum_{i,j} \omega_{i,j}^n  S(x-X_h(t^{n+1}; (x_i,v_j), t^n)S(v-V_h(t^{n+1};(x_i,v_j),t^n), \nonumber\\
& = & \sum_{i,j} \omega_{i,j}^{n+1}S(x-x_i)S(v-v_j),
\end{eqnarray}
by definition of $\omega_{i,j}^n$.
This implies a kind of continuity of $f_h$ at time $t^n$ on the grid points.

\subsubsection{Momentum conservation}

Let us precise that in this paragraph $(X_h(t;(x_i,v_j),t^n),V_h(t;(x_i,v_j),t^n))$ will be denoted $(X_{i,j}(t),V_{i,j}(t))$, and that Poisson will not be solved exactly.
The aim here is to prove that $\forall n$ 
\begin{eqnarray}\label{eq4.16}
\sum_{i,j}v_j f_h(t^n,x_i,v_j) = \sum_{i,j}v_j f_h(t^{n+1},x_i,v_j).
\end{eqnarray}

Let us distinguish two phases: the transport one and the deposition one. Let us start with the deposition phase, where we have to get:
 \begin{eqnarray*}
\sum_{i,j,k,l} \omega_{i,j}^{n+1}v_lS_1(x_k-x_i)S_1(v_l-v_j) = \sum_{i,j,k,l}\omega_{i,j}^n v_l S_1(x_k-X_{i,j}(t^{n+1}))S_1(v_l-V_{i,j}(t^{n+1}))
\end{eqnarray*}
\begin{eqnarray*}
\sum_{i,j,k,l} \omega_{i,j}^{n+1}v_lS_1(x_k-x_i)S_1(v_l-v_j) & = & \sum_{i,j} \omega_{i,j}^{n+1}\sum_l v_lS_1(v_l-v_j)\nonumber\\
							& = & \sum_{i,j} \omega_{i,j}^{n+1} v_j
\end{eqnarray*}
thanks to the  property (v) of linear splines. Moreover
\begin{eqnarray*}
\sum_l v_lS_1(v_l-V_{i,j}(t^{n+1}))  = V_{i,j}(t^{n+1}),
\end{eqnarray*}
thanks to the same property.
So we finally get for the deposition phase:
\begin{eqnarray}\label{eq4.17}
\sum_{i,j} \omega_{i,j}^{n+1} v_j=\sum_{i,j} \omega_{i,j}^nV(t^{n+1};(x_i,v_j),t^n).
\end{eqnarray}

\begin{remark}
This proof is given for linear splines, but was also checked for quadratic and cubic ones. The transport phase is independent from the spline degree.
\end{remark}

There remains to prove that 
\begin{eqnarray}\label{eq4.18}
\sum_{i,j} \omega_{i,j}^{n} v_j=\sum_{i,j} \omega_{i,j}^nV(t^{n+1};(x_i,v_j),t^n)
\end{eqnarray}
which corresponds to the transport phase. Note that this phase exists also in PIC methods, and the following proof of conservation of moments is adapted from \cite{birdsall}.

\paragraph{Verlet}

We have with our Verlet algorithm:
\begin{eqnarray}\label{eq4.19}
V(t^{n+1};(x_i,v_j),t^n)=v_j+\Delta t E^{n+\frac{1}{2}}(X(t^{n+\frac{1}{2}},(x_i,v_j),t^n)).
\end{eqnarray}
The electric field is only known on the mesh. In order to know it everywhere, we use a convolution between a spline function and the discrete $E$.
$$
E^{n+\frac{1}{2}}(X(t^{n+\frac{1}{2}},(x_i,v_j),t^n))=\sum_kE^{n+\frac{1}{2}}(x_k)S(x_k-X(t^{n+\frac{1}{2}},(x_i,v_j),t^n)).
$$
To get \eqref{eq4.18} using \eqref{eq4.19} we just have to prove that 
$$
\sum_{i,j,k}\omega_{i,j}^nE^{n+\frac{1}{2}}(x_k)S(x_k-X(t^{n+\frac{1}{2}},(x_i,v_j),t^n))=0,
$$
and 
\begin{eqnarray*}
\sum_{i,j,k}\omega_{i,j}^nE^{n+\frac{1}{2}}(x_k)S(x_k-X(t^{n+\frac{1}{2}},(x_i,v_j),t^n)) & = & \sum_kE^{n+\frac{1}{2}}(x_k)\sum_{i,j}\omega_{i,j}^nS(x_k-X(t^{n+\frac{1}{2}},(x_i,v_j),t^n)),\nonumber\\
		 & = &  \sum_kE^{n+\frac{1}{2}}(x_k)\rho^{n+\frac{1}{2}}(x_k),\nonumber\\
		 & = & 0,
\end{eqnarray*}
for most of the centered algorithms used to solve Poisson numerically, like the following centered finite difference one on staggered mesh, with linear regularization:
\begin{eqnarray*}
E^{n+\frac{1}{2}}(x_{k+\frac{1}{2}})-E^{n+\frac{1}{2}}(x_{k-\frac{1}{2}})=\Delta x \rho^{n+\frac{1}{2}}(x_{k}) \quad\forall k,
\end{eqnarray*}
and 
$$
E^{n+\frac{1}{2}}(x)=\sum_iE^{n+\frac{1}{2}}(x_{i})S_1(x-x_i).
$$
Indeed, we get:
\begin{eqnarray*}
\sum_kE^{n+\frac{1}{2}}(x_k)\rho^{n+\frac{1}{2}}(x_k) & = & \sum_k\sum_iE^{n+\frac{1}{2}}(x_k)E^{n+\frac{1}{2}}(x_i)((S_1(x_{k+\frac{1}{2}}-x_i)-S_1(x_{k+\frac{1}{2}}-x_i)),\nonumber\\
& = & \sum_iE^{n+\frac{1}{2}}(x_i)(\frac{E^{n+\frac{1}{2}}(x_{i-1})-E^{n+\frac{1}{2}}(x_{i+1})}{2}),\nonumber\\
&=& \frac{1}{2}(\sum_i E^{n+\frac{1}{2}}(x_i)E^{n+\frac{1}{2}}(x_{i+1})-\sum_i E^{n+\frac{1}{2}}(x_i)E^{n+\frac{1}{2}}(x_{i-1})),\nonumber\\
&=& 0,
\end{eqnarray*}
thanks to periodicity.

To conclude, using \eqref{eq4.17} and \eqref{eq4.18}, we get \eqref{eq4.16}, which is what was wanted.

\paragraph{CK algorithm.}
We still have to prove that:
\begin{eqnarray*}
\sum_{i,j} \omega_{i,j}^{n} v_j=\sum_{i,j} \omega_{i,j}^nV(t^{n+1};(x_i,v_j),t^n)
\end{eqnarray*}
That means for the third order scheme:
\begin{eqnarray*}
\sum_{i,j} \omega_{i,j}^{n}E^n(x_i) &= &0  \quad (i) \nonumber\\
\sum_{i,j} \omega_{i,j}^{n}\varphi^n(x_i,v_j) & = &0 \quad (ii) \nonumber\\
\sum_{i,j} \omega_{i,j}^{n}\phi^n(x_i,v_j) & = & 0 \quad (iii)\nonumber\\
\end{eqnarray*}
each number being linked with the order of the algorithm.

\paragraph{First order.}
Using the same strategy (regularization of the electric field and centered algorithm):
\begin{eqnarray}\label{eq4.20}
\sum_{i,j} \omega_{i,j}^{n}E^n(x_i) & = & \sum_{i,j,k,l} \omega_{i,j}^{n}E^n(x_k)S(x_i-x_k),\nonumber\\
	& = & \sum_{k,l}E^n(x_k)\rho^n(x_k) = 0.
\end{eqnarray}

\paragraph{Second order.}
For the second order, we need $J^n_i, \rho^n_i, \bar{J}$
$$
\rho^n_i=\Delta v \sum_{k,l}\omega_{k,l}^nS(x_i-x_k)-1,
$$
$$
J^n_i-\bar{J}=\Delta v\sum_{k,l}\omega_{k,l}^nv_lS(x_i-x_k)-\frac{\Delta x\Delta v}{L}\sum_{k,l}\omega_{k,l}^nv_l,
$$
\begin{eqnarray}\label{eq4.21}
\sum_{i,j}\omega_{i,j}^nv_j\rho^n(x_i)=\Delta v\sum_{i,j,k,l}\omega_{i,j}^n\omega_{k,l}^nv_jS(x_i-x_k)-\sum_{i,j}\omega_{i,j}^nv_j,
\end{eqnarray}
\begin{eqnarray}\label{eq4.22}
\sum_{i,j}\omega_{i,j}^nJ^n(x_i) & = & \Delta v\sum_{i,j,k,l}\omega_{i,j}^n\omega_{k,l}^nv_lS(x_i-x_k)-\frac{\Delta x\Delta v}{L}\sum_{i,jk,l}\omega_{k,l}^n\omega_{i,j}^nv_l,\nonumber\\
	& = &\Delta v\sum_{i,j,k,l}\omega_{i,j}^n\omega_{k,l}^nv_jS(x_i-x_k)-\sum_{k,l}\omega_{k,l}^nv_l,
\end{eqnarray}
using mass conservation and $\sum_{i,j}\omega_{i,j}^0=L$
\eqref{eq4.21} and \eqref{eq4.22} are the same, just exchanging $(i,j)$ and $(k,l)$.
So
$$
\sum_{i,j} \omega_{i,j}^{n}\varphi^n(x_i,v_j)  = 0.
$$

\paragraph{Third order.}
Here we need $I_2^n(x_i)$
\begin{eqnarray*}
I_2(t^n,x_i) & = & \int_{\mathbb{R}}v^2\sum_{k,l}\omega_{k,l}^nS(x_k-x_i)S(v-v_l)\,dv,\nonumber\\
	& = & \sum_{k,l}\omega_{k,l}^nS(x_k-x_i)\int_{\mathbb{R}}(v^2S(v)+2vv_lS(v)+ v_l^2S(v))\,dv,\nonumber\\
	& = & \sum_{k,l}\omega_{k,l}^nS(x_k-x_i)(\alpha +\Delta vv_l^2).
\end{eqnarray*}
In $\Phi$ we still have $\sum_{i,j} \omega_{i,j}^{n}E^n(x_i)  = 0$.
We also have three terms in $\frac{\partial}{\partial_x}$ which will be approached with a centered finite difference formula:
\begin{eqnarray}\label{eq4.23}
\sum_{i,j}\omega_{i,j}^n\frac{\partial}{\partial_x}I_2(t^n,x_i) & = & \frac{1}{2\Delta x}\sum_{i,j}\omega_{i,j}^n(I_2(t^n,x_{i+1})-I_2(t^n,x_{i-1})),\nonumber\\
& = &  \frac{1}{2\Delta x}\sum_{i,j,k,l}\omega_{i,j}^n\omega_{k,l}^n(\alpha+v_l^2)(S(x_{i+1}-x_k)-S(x_{i-1}-x_k)),
\end{eqnarray}
\begin{eqnarray}\label{eq4.24}
2\sum_{i,j}\omega_{i,j}^nv_j\frac{\partial J}{\partial_x}(t^n,x_i)  = \frac{1}{\Delta x}\sum_{i,j,k,l}\omega_{i,j}^n\omega_{k,l}^nv_jv_l(S(x_{i+1}-x_k)-S(x_{i-1}-x_k)),
\end{eqnarray}
\begin{eqnarray}\label{eq4.25}
\sum_{i,j}\omega_{i,j}^nv_j^2\frac{\partial \rho}{\partial_x}(t^n,x_i) =  \frac{1}{2\Delta x}\sum_{i,j,k,l}\omega_{i,j}^n\omega_{k,l}^nv_j^2(S(x_{i+1}-x_k)-S(x_{i-1}-x_k)).
\end{eqnarray}
Adding \eqref{eq4.23}, \eqref{eq4.24} and \eqref{eq4.25} and using \eqref{eq4.20} we have:
\begin{eqnarray*}
\sum_{i,j}\omega_{i,j}^n\phi^n(x_i)  =  \frac{1}{\Delta x}\sum_{i,j,k,l}\omega_{i,j}^n\omega_{k,l}^n(\alpha+v_l^2+2v_jv_l+v_j^2)(S(x_{i+1}-x_k)-S(x_{i-1}-x_k)),\nonumber\\
\end{eqnarray*}
\begin{eqnarray*}
\sum_{i,j,k,l}\omega_{i,j}^n\omega_{k,l}^n(\alpha+v_l^2+2v_jv_l+v_j^2)S(x_{i+1}-x_k) & = & \sum_{i,j,k,l}\omega_{k,l}^n\omega_{i,j}^n(\alpha+v_l^2+2v_jv_l+v_j^2)S(x_{k+1}-x_i)\nonumber\\
	 & =& \sum_{i,j,k,l}\omega_{i,j}^n\omega_{k,l}^n(\alpha+v_l^2+2v_jv_l+v_j^2)S(x_{k}-x_{i-1}), 
\end{eqnarray*}
just changing $(i,j)$ and $(k,l)$ and $S(x_{k+1}-x_i)=S(x_{k}-x_{i-1})$.
So we get:
\begin{eqnarray*}
\sum_{i,j}\omega_{i,j}^n\phi^n(x_i)=0.
\end{eqnarray*}


\begin{remark}
We can see that the conservation of the first moment in v implies that numerically $\frac{d\bar{J}}{dt}=0$, which means that $\frac{\Delta x \Delta v}{L}\sum_{i,j}^n\omega_{i,j}^nv_j$ is constant.
\end{remark}

\section{Convergence analysis}

\begin{theorem}
Assume that $f_0 \in W_{c,per_x}^{3,\infty}(\mathbb{R}_x \times \mathbb{R}_v)$, positive, periodic with respect to the variable $x$, with period $L$, and compactly supported in velocity.

Then the numerical solution of the Vlasov Poisson system $(f_h,E_h)$, computed by the numerical scheme introduced in section 3.2 converges towards the solution $(f,E)$ of the periodic Vlasov-Poisson system, and there exists a constant $C=C(\vert \vert f \vert \vert_{W^{1,\infty}(0,T;W^{2,\infty}(\Omega))})$ independent of $\Delta t$ and $h$ such that for Verlet and CK2 algorithms:
$$
\vert \vert f-f_h \vert \vert_{l^{\infty}(0,T;L^1(\Omega))} + \vert \vert E-E_h \vert \vert_{l^{\infty}(0,T;L^{\infty}([0,L]))} \leq C(\Delta t^2 + h^2 + \frac{h^2}{\Delta t}).
$$
For CK3,  we have:
$$
\vert \vert f-f_h \vert \vert_{l^{\infty}(0,T;L^1(\Omega))} + \vert \vert E-E_h \vert \vert_{l^{\infty}(0,T;L^{\infty}([0,L]))} \leq C(\Delta t^3 +h^2+ \frac{h^2}{\Delta t}).
$$
\begin{remark}
In order to get these estimates for CK, we will have to assume $\Delta t \leq \Delta x$.
\end{remark}
\end{theorem}

\subsection{Decomposition of the error}

Let $f$ be the exact solution of the Vlasov Poisson equation and $f_h$ the approximate solution previously defined.
In order to apply a discrete Gronwall inequality we express the $l^1$ error at time $t^{n+1}$
$$
e^{n+1}(i,j)=\vert f(t^{n+1},x_i,v_j)-f_h(t^{n+1},x_i,v_j) \vert \quad \forall (i,j),
$$
$$
e^{n+1}=\Delta x \Delta v \sum_{i,j}e^{n+1}(i,j).
$$
Then $f(t^{n+1},x_k,v_l)-f_h(t^{n+1},x_k,v_l)$ can be decomposed as 
\begin{eqnarray}\label{eq4.14}
 f(t^{n+1},x_k,v_l)-f_h(t^{n+1},x_k,v_l) & = &  f(t^{n+1},x_k,v_l)-R_hf(t^{n+1},x_k,v_l) + \nonumber\\
                                                                     &  &  R_hf(t^{n+1},x_k,v_l)-R_h\tilde{f_h}(t^{n+1},x_k,v_l) + \nonumber\\
                                                                      & &  R_h\tilde{f_h}(t^{n+1},x_k,v_l)-R_hf_h(t^{n+1},x_k,v_l),
\end{eqnarray}
where $\tilde{f_h}$ is the function $f_h$ at time $t^n$ but then  follows the exact characteristics. Since $f_h^{n+1}$ already belongs to the image of $R_h$, we have $R_hf_h(t^{n+1},x_k,v_l)=f_h(t^{n+1},x_k,v_l)$. 

In order to estimate $e^{n+1}$, the three terms of the right hand side of the previous equation have to be dealt with. These estimations are developed in the following subsection.

\subsection{A priori estimates}

\subsubsection{Stability for linear splines}

Let us translate the useful spline properties in this case, and give a few more results about the operator $R_h$.

\begin{lemma}
The $R_h$ operator is consistent, that is, using property (i), for $1\leq p \leq \infty$, and $0 \leq k \leq 1$
$$
\exists C > 0 \quad \vert \quad  \vert\vert f-R_hf \vert \vert_{W^{k,p}(\Omega)} \leq C h^{2-k} \vert f \vert_{W^{2,p}(\Omega)} \quad \forall f \in W^{2,p}(\Omega)\cap P(\Omega).
$$
\end{lemma}
This result is a classical property of B-splines.

\begin{lemma}
With linear splines, if $\omega_{i,j}(f_0) \geq 0  \quad \forall (i,j)$ \quad  then  $\forall n \quad \omega_{i,j}(f^n) \geq 0 \quad \forall (i,j)$. 
\end{lemma}

\noindent\textbf{Proof:} With linear interpolation, we get in fact $ \omega_{i,j}(f^n) = f^n(x_i,v_j) $, so if $f_0$ is positive, 
$\omega_{i,j}(f_0)$ is also, and since $ f^{n+1}(x_i,v_j)$ is a sum of positive contributions coming from the $f^n(x_k,v_l)$ which are positive by a recurrence hypothesis, it will also be positive, and so $\omega_{i,j}(f^{n+1})$ is positive for all $(i,j)$, and recurrently for all $n$.

\begin{lemma}
Stability: Let f belong to $C(\Omega) \cap P(\Omega)$, then we have:
$$
 \vert \vert R_hf \vert \vert_{L^1(\Omega)} = \vert \vert f \vert \vert_{L^1_h(\Omega)}.
$$
\end{lemma}

\noindent\textbf{Proof: }

\begin{eqnarray*}
\vert \vert R_hf \vert \vert_{L^1(\Omega)} & = & \int_0^L \int_{\mathbb{R}}  \vert R_hf(x,v) \vert  \,dv\,dx,\nonumber\\
							      & = &\int_0^L \int_{\mathbb{R}} \sum_{i,j} \omega_{i,j}(f)S(x-x_i)S(v-v_j),\nonumber\\
							      & = & \Delta x \Delta v \sum_{i,j} \omega_{i,j}(f),\nonumber\\
							      & = & \vert \vert f \vert \vert_{L^1_h(\Omega)},
\end{eqnarray*}
using $\int S(x)dx = \Delta x$, the positivity of $f$ thanks to Lemma 2 and the positivity of $f_0$.						      	 

\subsubsection{Towards Gronwall}

Let us precise that in this subsection, some lemmas are valid for all the time algorithms we use, and when they are not, the lemmas will be proved in each case successively. For the Cauchy Kovalevsky procedure, the proofs will be done for (CK 3), since their adaptation to lower orders is trivial.
We will now give estimates about the three right-hand side terms of the error $e^{n+1}$ \eqref{eq4.14}:

\begin{lemma} 
Let f belong to $C(\Omega) \cap P(\Omega)$, then we have:
\begin{eqnarray}\label{eq4.26}
\vert \vert f-R_hf \vert \vert_{L^1_h(\Omega)} \leq C h^2.
\end{eqnarray}

\end{lemma}

\noindent\textbf{Proof:}  Thanks to Lemma 3
\begin{eqnarray*}
\vert \vert f - R_hf \vert \vert_{L^1_h(\Omega)} &  = &  \vert \vert R_h(f-R_hf) \vert \vert_{L^1(\Omega)},\nonumber\\
								    & \leq & C  \vert \vert f-R_hf \vert \vert_{L^1(\Omega)},\nonumber\\	
								    & \leq & C'( \vert \vert f \vert \vert_{L^{\infty}(0,T;W^{2,\infty}(\Omega)})h^2,
\end{eqnarray*}
thanks to the property (ii) of spline interpolation and the fact that the domain is bounded.

\begin{lemma}
Let f belong to $C(\Omega) \cap P(\Omega)$, then we have:
\begin{eqnarray}\label{eq4.27}
\vert \vert R_hf^{n+1}-R_h\tilde{f_h}^{n+1} \vert \vert_{L^1_h(\Omega)} \leq e^n.
\end{eqnarray}
\end{lemma}

\noindent\textbf{Proof:} We compute
\begin{eqnarray*}
 \vert  \vert R_hf^{n+1} & - & R_h\tilde{f_h}^{n+1} \vert \vert_{L^1_h(\Omega)}  = \Delta x \Delta v \sum_{k,l} \vert (R_hf^{n+1}-R_h\tilde{f_h}^{n+1})(x_k,v_l) \vert ,\nonumber\\
			& = & \Delta x \Delta v \sum_{k,l} \vert \sum_{i,j} (f^{n+1}(x_i,v_j)-\tilde{f_h}^{n+1}(x_i,v_j))S(x_k-x_i)S(v_l-v_j) \vert, \nonumber\\
			& \leq & \Delta x \Delta v  \sum_{i,j} \vert \omega_{i,j}(f^n)-\omega_{i,j}(f_h^n) \vert \sum_{k,l} S(x_k-X(t^{n+1};(x_i,v_j),t^n)) S(v_l-V(t^{n+1};(x_i,v_j),t^n)), \nonumber\\
			& \leq & \Delta x \Delta v  \sum_{i,j} \vert \omega_{i,j}(f^n)-\omega_{i,j}(f_h^n) \vert,\nonumber\\
			& \leq & e^n,
\end{eqnarray*}													
thanks once more to the partition of unity (iii) and $f(x_i,v_j,t^n) =  \omega_{i,j}(f^n)$. 

\begin{lemma}
Let f belong to $C(\Omega) \cap P(\Omega)$, then we have:
\begin{eqnarray}\label{eq4.28}
\vert \vert R_h\tilde{f_h}^{n+1}-R_h{f_h}^{n+1} \vert \vert_{L^1_h(\Omega)} & \leq &  C \max_{i,j}(\vert X(t^{n+1};(x_i,v_j),t^n) - X_h(t^{n+1};(x_i,v_j),t^n) \vert ,\nonumber\\
& &+ \vert  V(t^{n+1};(x_i,v_j),t^n) - V_h(t^{n+1};(x_i,v_j),t^n) \vert).
\end{eqnarray}

\end{lemma}	

\noindent\textbf{Proof:}
\begin{eqnarray*}
& &\vert \vert R_h\tilde{f_h}^{n+1}  -  R_hf_h^{n+1} \vert \vert_{L^1_h(\Omega)}  =  \Delta x \Delta v \sum_{k,l} \vert (R_h\tilde{f_h}^{n+1}-R_h{f_h}^{n+1})(x_k,v_l) \vert ,\nonumber\\	
		& = & \Delta x \Delta v \sum_{k,l} \vert \sum_{i,j} \omega_{i,j}(f_h)^n( S(x_k-X(t^{n+1};(x_i,v_j),t^n)) S(v_l-V(t^{n+1};(x_i,v_j),t^n)),\nonumber\\
		& - & S(x_k-X_h(t^{n+1};(x_i,v_j),t^n)) S(v_l-V_h(t^{n+1};(x_i,v_j),t^n)))\vert .
\end{eqnarray*}
We can rewrite
\begin{multline*}
( S(x_k-X(t^{n+1};(x_i,v_j),t^n)) S(v_l-V(t^{n+1};(x_i,v_j),t^n))\nonumber\\
		- S(x_k-X_h(t^{n+1};(x_i,v_j),t^n)) S(v_l-V_h(t^{n+1};(x_i,v_j),t^n))) \nonumber\\
		 = ( S(x_k-X(t^{n+1};(x_i,v_j),t^n))-S(x_k-X_h(t^{n+1};(x_i,v_j),t^n)))S(v_l-V(t^{n+1};(x_i,v_j),t^n))\nonumber\\
		-(S(v_l-V(t^{n+1};(x_i,v_j),t^n))-S(v_l-V_h(t^{n+1};(x_i,v_j),t^n)))S(x_k-X_h(t^{n+1};(x_i,v_j),t^n)).
\end{multline*}	
Then, we use the fact that $S_1$  is 1-Lipschitzian, compactly supported, and the property (i):
\begin{multline*}
\sum_{k,l}\vert (S(x_k-X(t^{n+1};(x_i,v_j),t^n))-S(x_k-X_h(t^{n+1};(x_i,v_j),t^n)))S(v_l-V(t^{n+1};(x_i,v_j),t^n))\vert\nonumber\\
 \leq  \vert X(t^{n+1};(x_i,v_j),t^n) -  X_h(t^{n+1};(x_i,v_j),t^n) \vert ,
\end{multline*}
and
\begin{multline*}
\sum_{k,l}\vert (S(v_l-V(t^{n+1};(x_i,v_j),t^n))-S(v_l-V_h(t^{n+1};(x_i,v_j),t^n)))S(x_k-X_h(t^{n+1};(x_i,v_j),t^n))\vert\nonumber\\
\leq  \vert V(t^{n+1};(x_i,v_j),t^n) -  V_h(t^{n+1};(x_i,v_j),t^n) \vert .
\end{multline*}

So that we get:
\begin{multline*}
\vert \vert R_h\tilde{f_h}^{n+1}  -  R_hf_h^{n+1} \vert \vert_{L^1_h(\Omega)}		 \leq  \Delta x \Delta v \sum_{i,j} \vert  \omega_{i,j}(f_h)^n\vert (\vert X(t^{n+1};(x_i,v_j),t^n) -  X_h(t^{n+1};(x_i,v_j),t^n) \vert \nonumber\\
		 +  \vert  V(t^{n+1};(x_i,v_j),t^n) - V_h(t^{n+1};(x_i,v_j),t^n) )\vert),\nonumber\\
		 \leq  C \max_{i,j}(\vert X(t^{n+1};(x_i,v_j),t^n) - X_h(t^{n+1};(x_i,v_j),t^n) \vert + \vert  V(t^{n+1};(x_i,v_j),t^n) - V_h(t^{n+1};(x_i,v_j),t^n) \vert).
\end{multline*}			    	
 thanks to particle conservation ($\sum_{i,j}\omega_{i,j}(f_h^n)=\sum_{i,j}\omega_{i,j}(f_0))$ and positivity of $\omega_{i,j}(f_h^n)$, where $(X,V)$ are the exact characteristics, solution of the differential system \eqref{eq2.7}, and $(X_h,V_h)$ the approximate characteristics defined in \eqref{eq3.11}.
 
  To move on, we need another lemma which enables to control the difference between exact and computed characteristics. It clearly depends on the algorithm we use. Let us first give the lemma for the Verlet algorithm.
 
 \begin{lemma}: Verlet
 
 If $E \in W^{2,\infty}([0,t] \times \mathbb{R})$, and with $(X,V)$ calculated exactly with the differential system \eqref{eq2.7}, and $(X_h,V_h)$ computed with $E_h$ and a Verlet algorithm:
 \begin{eqnarray*}
 \vert X(t^{n+1};(x_i,v_j),t^n) & - & X_h(t^{n+1};(x_i,v_j),t^n) \vert +  \vert  V(t^{n+1};(x_i,v_j),t^n) - V_h(t^{n+1};(x_i,v_j),t^n) \vert\nonumber\\
  & \leq & C\Delta t^3  +  \Delta t\vert \vert(E-E_h)(t^{n+\frac{1}{2}}) \vert \vert_{L^{\infty}}.
\end{eqnarray*}

\end{lemma}

\noindent\textbf{Proof:} The strategy follows the work of M. Bostan and N. Crouseilles (\cite{bostan}).

Let us recall how $(x,v)^{n+1}$ is computed from $(x,v)^n$ with the Verlet algorithm
$$
x^{n+\frac{1}{2}}=x^n+\frac{\Delta t}{2}v^n,
$$
$$
v^{n+1}=v^n+\Delta t E_h(t^{n+\frac{1}{2}},x^{n}+\frac{\Delta t}{2}v^n),
$$
$$
x^{n+1}=x^{n+\frac{1}{2}}+\frac{\Delta t}{2}v^{n+1}.
$$
Then we define $X_h(t^{n+1},(x^n,v^n),t^n)=x^{n+1}$ and $V_h(t^{n+1},(x^n,v^n),t^n)=v^{n+1}$.

Let us begin with the characteristics in v:
\begin{eqnarray}\label{eq4.29}
(V_h-V)(t^{n+1};(x^n,v^n),t^n) & = &  v^n + \Delta t E_h(t^{n+\frac{1}{2}},x^{n}+\frac{\Delta t}{2}v^n) - v^n -\int_{t^n}^{t^{n+1}}E(s,X(s;(x^n,v^n),t^n)\,ds,\nonumber\\
   		& = & -\int_{t^n}^{t^{n+1}}(E(s,X(s;(x^n,v^n),t^n)-E(t^{n+\frac{1}{2}},X(t^{n+\frac{1}{2}};(x^n,v^n),t^n)))\,ds,\nonumber\\
		& - & \Delta t(E(t^{n+\frac{1}{2}},X(t^{n+\frac{1}{2}};(x^n,v^n),t^n))-E_h(t^{n+\frac{1}{2}},x^{n}+\frac{\Delta t}{2}v^n)).
\end{eqnarray}
Let us take care of the integral term in \eqref{eq4.29}, using a Taylor expansion around $s=t^{n+\frac{1}{2}}$ of $s\mapsto E(s,X(s;(x^n,v^n),t^n)$.
\begin{eqnarray}\label{eq4.30}
E(s,X(s;(x^n,v^n),t^n) & = & E(t^{n+\frac{1}{2}},X(t^{n+\frac{1}{2}};(x^n,v^n),t^n))+(s-t^{n+\frac{1}{2}})E'(t^{n+\frac{1}{2}},X(t^{n+\frac{1}{2}};(x^n,v^n),t^n))\nonumber\\
& + & \int_{t^{n+\frac{1}{2}}}^s(s-u)E''(u,X(u;(x^n,v^n),t^n)) \,du.
\end{eqnarray}
Let us precise that $E'(s,X(s))=\frac{d}{ds}E(s,X(s))$,  $E''(s,X(s))=\frac{d}{ds}E'(s,X(s))$. Then using \eqref{eq4.30} in \eqref{eq4.29}, we get
\begin{eqnarray*}
\int_{t^n}^{t^{n+1}}(E(s,X(s;(x^n,v^n),t^n))& - &E(t^{n+\frac{1}{2}},X(t^{n+\frac{1}{2}};(x^n,v^n),t^n))\,ds,\nonumber\\
&  = & E'(t^{n+\frac{1}{2}},X(t^{n+\frac{1}{2}};(x^n,v^n),t^n))[(\frac{s-t^{n+\frac{1}{2}}}{2})^2]_{t^n}^{t^{n+1}},\nonumber\\
					& + & \int_{t^n}^{t^{n+1}}\int_{t^{n+\frac{1}{2}}}^s(s-u)E''(u,X(u;(x,v),t^n)\, du\,ds.
\end{eqnarray*}
There are here two terms to control. The first one is zero, and for the second one, we have,
\begin{eqnarray}\label{eq4.31}
\vert  \int_{t^n}^{t^{n+1}}\int_{t^{n+\frac{1}{2}}}^s(s-u)E''(u,X(u;(x^n,v^n),t^n) \,du\,ds \vert & \leq & \| E'' \|_{L^{\infty}} \int_{t^n}^{t^{n+1}}\int_{t^{n+\frac{1}{2}}}^s(s-u)\,du\,ds ,\nonumber\\
			& \leq & \| E'' \|_{L^{\infty}} \int_{t^n}^{t^{n+1}}[-(\frac{s-u}{2})^2]_{t^{n+\frac{1}{2}}}^s \,ds,\nonumber\\
			& \leq & \frac{1}{2}\| E'' \|_{L^{\infty}}\int_{t^n}^{t^{n+1}}(s-t^{n+\frac{1}{2}})^2\,ds, \nonumber\\
			& \leq & \frac{\Delta t^3}{24}\| E'' \|_{L^{\infty}} \leq C\| E'' \|_{L^{\infty}}\Delta t^3.
\end{eqnarray}
Now let us deal with the second term of \eqref{eq4.29}. Since E is bounded:
\begin{eqnarray}\label{eq4.32}
\vert x^n+ \frac{\Delta t}{2}v^n-X(t^{n+\frac{1}{2}};(x^n,v^n),t^n) \vert & = & \vert \int_{t^n}^{t^{n+\frac{1}{2}}} v^n - V(s;(x^n,v^n),t^n)\,ds \vert,\nonumber\\
			& \leq &  \int_{t^n}^{t^{n+\frac{1}{2}}}(s-t^n)V'(u;(x^n,v^n),t^n)\,ds \quad u\in[t^n,s],\nonumber\\
			& \leq & C(\vert\vert E \vert\vert_{L^{\infty}}) \int_{t^n}^{t^{n+\frac{1}{2}}} (s-t^n)\,ds \leq C'\Delta t^2.
\end{eqnarray}
and thus, since $E'$ is bounded, using the zero mean theorem:
\begin{eqnarray}\label{eq4.33}
 \vert E(t^{n+\frac{1}{2}},X(t^{n+\frac{1}{2}};(x^n,v^n),t^n)) - E(t^{n+\frac{1}{2}},x^{n+\frac{1}{2}}) \vert & \leq & C(\vert\vert E' \vert\vert_{L^{\infty}}) \vert X(t^{n+\frac{1}{2}};(x^n,v^n),t^n)-x^{n+\frac{1}{2}} \vert ,\nonumber\\
 		& \leq & C(\vert\vert E' \vert\vert_{L^{\infty}})\Delta t^2.
\end{eqnarray}
Finally the second term of \eqref{eq4.29} can be controlled by
\begin{eqnarray}\label{eq4.34}
\vert E(t^{n+\frac{1}{2}},X(t^{n+\frac{1}{2}};(x^n,v^n),t^n))& - &E_h(t^{n+\frac{1}{2}},x^{n+\frac{1}{2}}) \vert  \leq  \vert \vert (E-E_h)(t^{n+\frac{1}{2}}) \vert \vert_{\infty}\nonumber\\ 
		& + & \vert E(t^{n+\frac{1}{2}},X(t^{n+\frac{1}{2}};(x^n,v^n),t^n))-E(t^{n+\frac{1}{2}},x^{n+\frac{1}{2}}) \vert.
\end{eqnarray}
So, using \eqref{eq4.31}, \eqref{eq4.33} and \eqref{eq4.34}, we get
\begin{eqnarray}\label{eq4.35}
\vert (V_h-V)(t^{n+1};(x,v),t^n) \vert \leq C\Delta t ^3 + \Delta t \vert \vert (E-E_h)(t^{n+\frac{1}{2}}) \vert \vert_{L^{\infty}(\Omega)}.
\end{eqnarray}
Let us now deal with the characteristics in $X$:
\begin{eqnarray}\label{eq4.36}
(X_h-X)(t^{n+1};(x^n,v^n),t^n) & = & \frac{\Delta t}{2}v^n + \frac{\Delta t}{2}v^{n+1} - \int_{t^n}^{t^{n+1}}V(s;(x^n,v^n),t^n)\,ds,\nonumber\\
		& = & - \Delta t(V(t^{n+\frac{1}{2}};(x^n,v^n),t^n)-\frac{1}{2}v^n-\frac{1}{2}v^{n+1}),\nonumber\\
		& - & \int_{t^n}^{t^{n+1}}(V(s;(x^n,v^n),t^n)-V(t^{n+\frac{1}{2}};(x^n,v^n),t^n)))\,ds,
\end{eqnarray}
so once again we have to control two terms.

For the first one, thanks to Taylor's inequality, like for X, it comes:
\begin{eqnarray}\label{eq4.37}
\vert \int_{t^n}^{t^{n+1}}(V(s;(x^n,v^n),t^n)-V(t^{n+\frac{1}{2}};(x^n,v^n),t^n))\,ds \vert \leq C(\vert\vert E \vert\vert_{L^{\infty}})\Delta t^3.
\end{eqnarray}
Let us precise that this is nothing else than the error in  the mid-point rule for numerical integration.
Now, the second term in \eqref{eq4.36}:
\begin{eqnarray*}
V(t^{n+\frac{1}{2}};(x^n,v^n),t^n)) & = & v^n + \int_{t^n}^{t^{n+\frac{1}{2}}}E(s,X(s;(x^n,v^n),t^n))\,ds,\nonumber\\
					   & = & v^n+ \frac{\Delta t}{2} E(t^{n+\frac{1}{2}},X(t^{n+\frac{1}{2}};(x^n,v^n),t^n)),\nonumber\\
					   & + &  \int_{t^n}^{t^{n+\frac{1}{2}}}(E(s,X(s;(x^n,v^n),t^n))- E(t^{n+\frac{1}{2}},X(t^{n+\frac{1}{2}};(x^n,v^n),t^n)))\,ds,
\end{eqnarray*}
and					   
$$
\vert \int_{t^n}^{t^{n+\frac{1}{2}}}(E(s,X(s;(x^n,v^n),t^n))- E(t^{n+\frac{1}{2}},X(t^{n+\frac{1}{2}};(x^n,v^n),t^n)))\,ds \vert \leq C(\vert\vert E' \vert\vert_{L^{\infty}})\Delta t^2, 
 $$
with the error formula for the rectangle rule. On the other hand
 \begin{eqnarray*}
 V(t^{n+\frac{1}{2}};(x^n,v^n),t^n)-\frac{1}{2}v^n-\frac{1}{2}v^{n+1} & = & \frac{1}{2}v^n-\frac{1}{2}v^{n+1}\nonumber\\
   		& + & \frac{\Delta t}{2} E(t^{n+\frac{1}{2}},X(t^{n+\frac{1}{2}};(x^n,v^n),t^n)) +O(\Delta t^2).
\end{eqnarray*}
Since 
$$
\frac{1}{2}(v^n-v^{n+1})= -\frac{\Delta t}{2} E_h(t^{n+\frac{1}{2}},x^{n+\frac{1}{2}}),
$$
we have, proceeding as for \eqref{eq4.33}-\eqref{eq4.34}
\begin{eqnarray}\label{eq4.38}
\vert  V(t^{n+\frac{1}{2}};(x^n,v^n),t^n)-\frac{1}{2}v^n-\frac{1}{2}v^{n+1} \vert & = & \vert \frac{\Delta t}{2}(E(t^{n+\frac{1}{2}},X(t^{n+\frac{1}{2}};(x^n,v^n),t^n))-E_h(t^{n+\frac{1}{2}},x^{n+\frac{1}{2}})\nonumber\\
&& \hspace{5cm}+O(\Delta t^2)\vert,\nonumber\\
		& \leq & \frac{\Delta t}{2}(\vert \vert (E-E_h)(t^{n+\frac{1}{2}}) \vert \vert_{L^{\infty}(\Omega)}+C\Delta t^2).
\end{eqnarray}

To conclude, using \eqref{eq4.37}, \eqref{eq4.38}, we have:
\begin{eqnarray}\label{eq4.39}
\vert (X_h-X)(t^{n+1};(x,v),t^n) \vert \leq C\Delta t ^3 + C \Delta t^2 \vert \vert (E-E_h)(t^{n+\frac{1}{2}}) \vert \vert_{L^{\infty}(\Omega)}.
\end{eqnarray}
Finally, using \eqref{eq4.35} and \eqref{eq4.39}, we get the estimation of Lemma 7, and using \eqref{eq4.28}, this also implies
 \begin{eqnarray*}
\vert \vert R_h\tilde{f_h}^{n+1}  -  R_hf_h^{n+1} \vert \vert_{L^1_h(\Omega)} \leq C\Delta t^3  + \Delta t\vert \vert(E-E_h)(t^{n+\frac{1}{2}}) \vert \vert_{L^{\infty}}.
\end{eqnarray*}

\begin{lemma}: CK3

If $E \in W^{4,\infty}([0,t] \times \mathbb{R})$, and with $(X,V)$ calculated exactly with the differential system \eqref{eq2.7}, and $(X_h,V_h)$ computed with $E_h, \rho_h, J_h$ and a CK3 algorithm:
 \begin{eqnarray*}
 \vert X(t^{n+1};(x_i,v_j),t^n) & - & X_h(t^{n+1};(x_i,v_j),t^n) \vert +  \vert  V(t^{n+1};(x_i,v_j),t^n) - V_h(t^{n+1};(x_i,v_j),t^n) \vert\nonumber\\
  & \leq & C\Delta t^4 + C(\Delta t\vert \vert(E^n-E_h^n)\vert \vert_{l^{\infty}(\Omega}  + \Delta t^2\vert \vert(\phi^n-\phi^n_h) \vert \vert_{l^{\infty}(\Omega)}\nonumber\\
  & & +\Delta t^3\vert \vert(\varphi^n-\varphi^n_h) \vert \vert_{l^{\infty}(\Omega)}).
\end{eqnarray*}
\end{lemma}

\noindent\textbf{Proof:}
This proof just relies on Taylor expansions and computations already made:
$$
X(t^{n+1};(x_i,v_j),t^n)=x_i+\Delta t v_j + \frac{\Delta t ^2}{2}E^n(x_i) +  \frac{\Delta t ^3}{6}\phi^n(x_i,v_j) + O(\Delta t ^4),
$$
$$
X_h(t^{n+1};(x_i,v_j),t^n)=x_i+\Delta t v_j + \frac{\Delta t ^2}{2}E_h^n(x_i) +  \frac{\Delta t ^3}{6}\phi_h^n(x_i,v_j),
$$
and
$$
V(t^{n+1};(x_i,v_j),t^n)=v_j + \Delta t E^n(x_i) +  \frac{\Delta t ^2}{2} \phi^n(x_i,v_j)  + \frac{\Delta t ^3}{6}\varphi^n(x_i,v_j) + O(\Delta t ^4),
$$
$$
V_h(t^{n+1};(x_i,v_j),t^n)=v_j + \Delta t E_h^n(x_i) +  \frac{\Delta t ^2}{2} \phi_h^n(x_i,v_j)  + \frac{\Delta t ^3}{6}\varphi_h^n(x_i,v_j),
$$
and the lemma follows by simple subtraction.

In both cases we need to control the difference between the exact and approximate fields. Let us begin with Verlet algorithm.

\begin{lemma}: Verlet

If $E \in W^{2,\infty}([0,t] \times \mathbb{R})$ , it comes
$$
 \vert \vert (E-E_h)(t^{n+\frac{1}{2}}) \vert \vert_{L^{\infty}(\Omega)} \leq C(h^2 + \Delta t^2 +\Delta t \, h^2 + e^n) 
$$

\end{lemma}

\noindent\textbf{Proof:}
First
$$
E(t^{n+\frac{1}{2}},x)= \int_0^LK(x,y)(\int_{\mathbb{R}}f(t^{n+\frac{1}{2}},y,v)dv-1)\,dy,
$$
$$
E_h(t^{n+\frac{1}{2}},x)= \int_0^LK(x,y)(\int_{\mathbb{R}}f_h(t^{n+\frac{1}{2}},y,v)dv-1)\,dy.
$$
Hence
\begin{eqnarray}\label{eq4.40}
E(t^{n+\frac{1}{2}},x)-E_h(t^{n+\frac{1}{2}},x) & = & \int_0^LK(x,y)(\int_{\mathbb{R}}(f(t^{n+\frac{1}{2}},y,v)-f_h(t^{n+\frac{1}{2}},y,v))\,dv)\,dy,\nonumber\\
				& = &  \int_0^LK(x,y)(\int_{\mathbb{R}}(f(t^{n+\frac{1}{2}},y,v)-R_hf(t^{n+\frac{1}{2}},y,v))\,dv)\,dy,\nonumber\\
				& + &  \int_0^LK(x,y)(\int_{\mathbb{R}}(R_hf(t^{n+\frac{1}{2}},y,v)-\bar{f}(t^{n+\frac{1}{2}},y,v))\,dv)\,dy,\nonumber\\
				& + & \int_0^LK(x,y)(\int_{\mathbb{R}}(\bar{f}(t^{n+\frac{1}{2}},y,v)-\tilde{f_h}(t^{n+\frac{1}{2}},y,v))\,dv)\,dy,\nonumber\\
				& + & \int_0^LK(x,y)(\int_{\mathbb{R}}(\tilde{f_h}(t^{n+\frac{1}{2}},y,v)-f_h(t^{n+\frac{1}{2}},y,v))\,dv)\,dy,
\end{eqnarray}
where 
$$\bar{f}(t^{n+\frac{1}{2}},y,v)=\sum_{k,l}\omega_{k,l}(f^n)S(y-X(t^{n+\frac{1}{2}};(x_k,v_l),t^n)S(v-V(t^{n+\frac{1}{2}};(x_k,v_l),t^n)$$ and $$\tilde{f_h}(t^{n+\frac{1}{2}},y,v)=\sum_{k,l}\omega_{k,l}(f_h^n)S(y-X(t^{n+\frac{1}{2}};(x_k,v_l),t^n)S(v-V(t^{n+\frac{1}{2}};(x_k,v_l),t^n).$$ 
In order to lighten notations, $X(t^{n+\frac{1}{2}};(x_k,v_l),t^n)$ and $V(t^{n+\frac{1}{2}};(x_k,v_l),t^n)$ will be denoted $X_{k,l}^{n+\frac{1}{2}}$ and $V_{k,l}^{n+\frac{1}{2}}$. We have four terms to control.

The first one is controlled using property (ii) of consistency and accuracy:
\begin{eqnarray}\label{eq4.41}
\vert  \int_0^LK(x,y)(\int_{\mathbb{R}}(f(t^{n+\frac{1}{2}},y,v)-R_hf(t^{n+\frac{1}{2}},y,v))\,dv)\,dy \vert  \leq  C \vert\vert K \vert\vert_{\infty}h^2.
\end{eqnarray}

Now, the second term of \eqref{eq4.40}. 
\begin{multline*}
\int_0^LK(x,y)(\int_{\mathbb{R}}(R_hf(t^{n+\frac{1}{2}},y,v)-\bar{f}(t^{n+\frac{1}{2}},y,v))\,dv)\,dy,\nonumber\\
=\int_0^LK(x,y)(\int_{\mathbb{R}}(\sum_{k,l}(\omega_{k,l}^{n+\frac{1}{2}}S(y-x_k)S(v-v_l)-\omega_{k,l}^nS(y-X_{k,l}^{n+\frac{1}{2}})S(v-V_{k,l}^{n+\frac{1}{2}})))\,dv)\,dy.
\end{multline*}
We have, using Taylor expansion and Vlasov equation:
\begin{eqnarray*}
\omega_{k,l}^{n+\frac{1}{2}} & = & f^{n+\frac{1}{2}}(x_k,v_l)=f^n(x_k,v_l)+\frac{\Delta t}{2}\frac{\partial f}{\partial t}(x_k,v_l)+O(\Delta t^2),\nonumber\\
& = & f^n(x_k,v_l)-\frac{\Delta t}{2}(v_l \frac{\partial f}{\partial x}(x_k,v_l)+E^n(x_k) \frac{\partial f}{\partial v}(x_k,v_l))+O(\Delta t^2).
\end{eqnarray*}
Moreover, since S is piecewise polynomial of degree one and continuous, we have almost everywhere (which is enough as we are going to integrate these expressions)
$$
S(y-X_{k,l}^{n+\frac{1}{2}})=S(y-x_k-\frac{\Delta t}{2}v_l+O(\Delta t^2))=S(y-x_k)-\frac{\Delta t}{2}v_lS'(y-x_k)+O(\Delta t^2).
$$
$$
S(v-V_{k,l}^{n+\frac{1}{2}})=S(v-v_l-\frac{\Delta t}{2}E^n(x_k)+O(\Delta t^2))=S(v-v_l)-\frac{\Delta t}{2}v_lS'(v-v_l)+O(\Delta t^2).
$$
Using, these expansions, we get:
\begin{align}\label{term2}
&\int_0^LK(x,y)(\int_{\mathbb{R}}(\sum_{k,l}(\omega_{k,l}^{n+\frac{1}{2}}S(y-x_k)S(v-v_l)-\omega_{k,l}^nS(y-X_{k,l}^{n+\frac{1}{2}})S(v-V_{k,l}^{n+\frac{1}{2}})))\,dv)\,dy,\nonumber\\
 &=  \frac{\Delta t}{2}\big( \int_0^LK(x,y)(\int_{\mathbb{R}}(\sum_{k,l} v_l S(v-v_l)(f^n(x_k,v_l)S'(y-x_k)-\frac{\partial f^n}{\partial x}(x_k,v_l)S(y-x_k)),\nonumber\\
 &+   \sum_{k,l} E^n(x_k)S(y-x_k)(f^n(x_k,v_l)S'(v-v_l)-\frac{\partial f^n}{\partial x}(x_k,v_l)S(v-v_l)))\,dv)\,dy\big) + O(\Delta t^2).
\end{align}

There are two terms in \eqref{term2}. They will be dealt with similarly using mid-point quadrature, which is of second order. For the first term, it writes
\begin{align*}
&  \int_0^LK(x,y)(\int_{\mathbb{R}}(\sum_{k,l} v_l S(v-v_l)(f^n(x_k,v_l)S'(y-x_k)-\frac{\partial f^n}{\partial x}(x_k,v_l)S(y-x_k))\,dv\,dy,\nonumber\\
& =  \Delta v \sum_{k,l}v_l\int_0^LK(x,y)(f^n(x_k,v_l)S'(y-x_k)-\frac{\partial f^n}{\partial x}(x_k,v_l)S(y-x_k))\,dy,\nonumber\\
& =   \Delta v \sum_{i,k,l}(v_lK(x,x_{i+\frac{1}{2}})(f^n(x_k,v_l)S'(x_{i+\frac{1}{2}}-x_k)-\frac{\partial f^n}{\partial x}(x_k,v_l)S(x_{i+\frac{1}{2}}-x_k)).
\end{align*}
Here, we have to use the properties of linear splines. $S'(x_{i+\frac{1}{2}}-x_k)$ and $S(x_{i+\frac{1}{2}}-x_k)$ are non zero only if $k=i$ or $k=i+1$. Then, we have: $S'(x_{i+\frac{1}{2}}-x_i)=\frac{-1}{\Delta x}$, $S'(x_{i+\frac{1}{2}}-x_{i+1})=\frac{1}{\Delta x}$ and $S(x_{i+\frac{1}{2}}-x_i)=S(x_{i+\frac{1}{2}}-x_{i+1})=\frac{1}{2}$. Using that, we get 
\begin{align*}
&  \Delta v \sum_{i,k,l}(v_lK(x,x_{i+\frac{1}{2}})(f^n(x_k,v_l)S'(x_{i+\frac{1}{2}}-x_k)-\frac{\partial f^n}{\partial x}(x_k,v_l)S(x_{i+\frac{1}{2}}-x_k)),\nonumber\\
& =   \Delta v \sum_{i,l}v_lK(x,x_{i+\frac{1}{2}})\left( \frac{f^n(x_{i+1},v_l)-f^n(x_i,v_l)}{\Delta x}-\frac{1}{2}\left(\frac{\partial f^n}{\partial x}(x_i,v_l)+\frac{\partial f^n}{\partial x}(x_{i+1},v_l)\right)\right).
\end{align*}
Using Taylor expansions with respect to x, we easily get:
$$
\frac{f(x_{i+1},v_l)-f(x_i,v_l)}{\Delta x}=\frac{\partial f^n}{\partial x}(x_{i+\frac{1}{2}},v_l)+O(\Delta x^2).
$$
and
$$
\frac 12\left(\frac{\partial f^n}{\partial x}(x_i,v_l)+\frac{\partial f^n}{\partial x}(x_{i+1},v_l)\right)=\frac{\partial f^n}{\partial x}(x_{i+\frac{1}{2}},v_l)+O(\Delta x^2),
$$
since $f \in W_{c,per_x}^{3,\infty}(\mathbb{R}_x \times \mathbb{R}_v)$.
And to conclude for the first term of \eqref{term2}:
$$
\vert \int_0^LK(x,y)(\int_{\mathbb{R}}(\sum_{k,l} v_l S(v-v_l)(f^n(x_k,v_l)S'(y-x_k)-\frac{\partial f^n}{\partial x}(x_k,v_l)S(y-x_k))\,dv\,dy\vert \leq C \Delta x^2.
$$
For the second term of  \eqref{term2}, with a mid-point quadrature for the integral with respect to v, and the same properties of splines:
\begin{align*}
& \int_0^LK(x,y)(\int_{\mathbb{R}}(\sum_{k,l} E^n(x_k)S(y-x_k)(f^n(x_k,v_l)S'(v-v_l)-\frac{\partial f^n}{\partial x}(x_k,v_l)S(v-v_l)))\,dv)\,dy,\nonumber\\
& =  \Delta v\sum_{k,l}E^n(x_k)\int_0^LK(x,y)S(y-x_k)\,dy \sum_j(f^n(x_k,v_l)S'(v_{j+\frac{1}{2}}-v_l)-\frac{\partial f^n}{\partial x}(x_k,v_l)S(v_{j+\frac{1}{2}}-v_l)),\nonumber\\
& \leq \Delta x\,\Delta v \|K\|_{L^\infty}  \sum_{k,j}|E^n(x_k)|\big|\frac{f^n(x_k,v_{j+1})-f^n(x_k,v_j)}{\Delta v}-\frac{1}{2}(\frac{\partial f^n}{\partial v}(x_k,v_{j+1})+\frac{\partial f^n}{\partial v}(x_k,v_j))\big|.
\end{align*}
Using again Taylors expansions, with respect to v, we get:
\begin{align*}
& \vert \int_0^LK(x,y)(\int_{\mathbb{R}}(\sum_{k,l} E^n(x_k)S(y-x_k)(f^n(x_k,v_l)S'(v-v_l)-\frac{\partial f^n}{\partial x}(x_k,v_l)S(v-v_l)))\,dv)\,dy,\vert \nonumber\\
& \leq  \Delta x \Delta v \vert \vert K \vert \vert_{L^{\infty}}\sum_{k,j}\vert E^n(x_k) \vert \Delta v^2,\nonumber\\
& \leq  C \vert \vert E^n \vert \vert_{L_\infty([0,L])} \Delta v^2 \leq C' \Delta v^2,
\end{align*}
since $E \in W^{2,\infty}([0,t] \times \mathbb{R})$.
To conclude, the second term of \eqref{eq4.40} can be bounded like that:
\begin{eqnarray}\label{eq4.42}
 \vert \int_0^LK(x,y)(\int_{\mathbb{R}}(R_hf(t^{n+\frac{1}{2}},y,v)-\bar{f}(t^{n+\frac{1}{2}},y,v))\,dv)\,dy\vert \leq C(\Delta t \,h^2 +\Delta t^2).
\end{eqnarray}

For the third term of \eqref{eq4.40}:
\begin{align}\label{eq4.43}
&  \vert \int_0^LK(x,y)(\int_{\mathbb{R}}(\bar{f}(t^{n+\frac{1}{2}},y,v)-\tilde{f_h}(t^{n+\frac{1}{2}},y,v))\,dv)\,dy, \vert \nonumber\\
& =  \vert \int_0^LK(x,y)(\int_{\mathbb{R}} \sum_{k,l} (\omega_{k,l}(f^n)-\omega_{k,l}(f_h^n))S(y-X_{k,l}^{n+\frac{1}{2}})S(v-V_{k,l}^{n+\frac{1}{2}})))\,dv)\,dy,\nonumber\\
& \leq  \Delta v \sum_{k,l} \vert \omega_{k,l}(f^n)-\omega_{k,l}(f_h^n) \vert \int_0^L \vert K(x,y) \vert S(y-X_{k,l}^{n+\frac{1}{2}})\,dy,\nonumber\\
&\leq  \vert \vert K \vert \vert_{\infty} \Delta x \Delta v \sum_{k,l} \vert f^n(x_k,v_l)-f_h^n(x_k,v_l) \vert \leq C e^n.
\end{align}

Eventually, the last term of \eqref{eq4.40}:
\begin{align}\label{eq4.44}
&  \vert \int_0^LK(x,y)(\int_{\mathbb{R}}(\tilde{f_h}(t^{n+\frac{1}{2}},y,v)-f_h(t^{n+\frac{1}{2}},y,v))\,dv)\,dy \vert, \nonumber\\
& =  \vert  \int_0^LK(x,y)(\int_{\mathbb{R}}\sum_{k,l}\omega_{k,l}(f_h^n)\big( S(y-X_{k,l)}^{n+\frac{1}{2}})S(v-V_{k,l}^{n+\frac{1}{2}})-S(y-X_{h;k,l}^{n+\frac{1}{2}})S(v-V_{h;k,l}^{n+\frac{1}{2}})\big) \,dv)\,dy \vert, \nonumber\\
& \leq  \Delta v \sum_{k,l}\omega_{k,l}(f_h^n) \int_0^L K(x,y)\big(S(y-X_{k,l}^{n+\frac{1}{2}})-S(y-X_{h;k,l)}^{n+\frac{1}{2}})\vert\,dy, \nonumber\\
& \leq  \vert \vert K \vert \vert_{\infty} \Delta x \Delta v \sum_{k,l}\omega_{k,l}(f_h^n) \vert X_{k,l}^{n+\frac{1}{2}}-X_{h;k,l}^{n+\frac{1}{2}}\vert, \nonumber\\
& \leq  C \Delta t^2,
\end{align}
using successively the positivity of $f_h$, mass conservation, $f_h$ is 1-lipschitzian and a result in Lemma 7: \eqref{eq4.32}, where obviously, $(X_{h;(k,l)}^{n+\frac{1}{2}},V_{h;(k,l)}^{n+\frac{1}{2}})$ are the appoximate characteristic curves at time $t^{n+\frac{1}{2}}$ beginning at time $t^n$ at $(x_k,v_l)$.

To conclude, using \eqref{eq4.41}, \eqref{eq4.42}, \eqref{eq4.43} and \eqref{eq4.44} we get:
$$
 \vert \vert (E-E_h)(t^{n+\frac{1}{2}}) \vert \vert_{L^{\infty}(\Omega)} \leq C(h^2 + \Delta t^2 +e^n + \Delta t \,h^2),
$$
which is what was expected.
 
\begin{lemma} CK3
With the same hypothesis as in Lemma 8, we have :

\begin{eqnarray*}
\vert \vert R_h\tilde{f_h}^{n+1}-R_h{f_h}^{n+1} \vert \vert_{L^1_h(\Omega)} & \leq &  C(e^n(\Delta t + \frac{\Delta t^3}{\Delta x^2}+\frac{\Delta t^2}{\Delta x} ) + \Delta t\, h^2 + \Delta t^2 \frac{h^2}{\Delta x} + \Delta t^3 \frac{h^2}{\Delta x^2} + \Delta t^4).
\end{eqnarray*}
\end{lemma}

\noindent\textbf{Proof:}
Here, we need to evaluate the difference between the $l^{\infty}(\Omega)$ norms of the exact and approximate values of $\phi,\varphi $ and $\psi$, so to say the one between:
\begin{itemize}
\item  $E^n$ and $E_h^n$, $\rho^n$ and $\rho_h^n$, $J^n$ and $J_h^n$,
\item  their first spatial derivative and the one of $I_2^n$ and $I_{2  h}^n$.
\end{itemize}
Let us start with
\begin{eqnarray*}
(E^n-E_h^n)(x_i) & = &  \int_0^LK(x_i,y)(\int_{\mathbb{R}}(f(t^n,y,v)-f_h(t^n,y,v)dv)\,dy,\nonumber\\
			    & \leq & C(e^n +h^2).
\end{eqnarray*}
simply using a quadrature with the mesh points, which will also be of second order thanks to periodicity, and the fact that $K$ is bounded.
So that:
\begin{eqnarray}\label{eq4.46}
\vert \vert E^n-E_h^n \vert \vert_{l^{\infty}([0,L])} \leq C (e^n + h^2),
\end{eqnarray}
\begin{eqnarray*}
 \vert \vert \rho^n-\rho_h^n \vert \vert_{L_h^1([0,L])} & = & \Delta x \sum_i \vert \int_{\mathbb{R}}(f^n-f_h^n)(x_i,v)dv\vert,\nonumber\\
 		& \leq & C(e^n + h^2).
\end{eqnarray*}
So that using the equivalence of discrete norms, carefully noticing that $$\vert \vert . \vert \vert_{l^{\infty}([0,L])} \leq \frac{1}{\Delta x}\vert \vert . \vert \vert_{L_h^{1}([0,L])}$$  we get:
\begin{eqnarray*}
 \vert \vert \rho^n-\rho_h^n \vert \vert_{l^{\infty}([0,L])} & \leq &  \frac{1}{\Delta x}  \vert \vert \rho^n-\rho_h^n \vert \vert_{L_h^1([0,L])},\nonumber\\
 & \leq & \frac{C}{\Delta x}(e^n + h^2),
\end{eqnarray*}
and
\begin{eqnarray*}
 \vert \vert J^n-J_h^n \vert \vert_{L_h^1([0,L])} & = & \Delta x \sum_i \vert \int_{\mathbb{R}}v(f^n-f_h^n)(x_i,v)dv\vert,\nonumber\\
 		& \leq & C(e^n + h^2),
\end{eqnarray*}
using the same arguments and the fact that f is compactly supported, so that: 
\begin{eqnarray*}
 \vert \vert J^n-J_h^n \vert \vert_{l^{\infty}([0,L])} & \leq & \frac{1}{\Delta x} \vert \vert J^n-J_h^n \vert \vert_{L_h^1([0,L])},\nonumber\\
 & \leq & \frac{C}{\Delta x}(e^n + h^2).
 \end{eqnarray*}
 Let us precise that the same bound is obviously also valid for $\bar{J}$.
 So that we get, still using that $\Omega$ is bounded:
 $$
\vert \vert \phi^n-\phi_h^n \vert \vert_{l^{\infty}([0,L])}\leq  \frac{C}{\Delta x}(e^n + h^2).
$$

 For the first spatial derivative of these three terms, we can use the same strategy of finite difference. Let us do it with $E$.
 \begin{eqnarray*}
 \frac{\partial (E-E_h)}{\partial x}(t^n,x_i)= \frac{(E^n_{i+1}-E^n_{i+1, h})-(E^n_{i-1}-E^n_{i-1, h})}{2\Delta x} + O(\Delta x^2),
 \end{eqnarray*}
 so that using \eqref{eq4.46} we get:
 \begin{eqnarray*}
 \vert \vert \frac{\partial (E-E_h)}{\partial x}(t^n) \vert \vert_{l^{\infty}([0,L])}  \leq  \frac{C}{\Delta x}(e^n + h^2).
 \end{eqnarray*}
For $\rho^n,J^n$ and for $I_2^n$ just bounding $v,v^2$ in its integral definition, the same strategy leads to
 $$
\vert \vert \varphi^n-\varphi_h^n \vert \vert_{l^{\infty}([0,L])}\leq \frac{C}{\Delta x^2}(e^n + h^2).
$$
Plugging these estimates into Lemma 8 and then Lemma 6 completes the proof.
 
\subsection{End of the proof}

For the sake of simplicity, and since we are interested in $\Delta t, \Delta x $ tend to $0$, we will assume $\Delta t \leq 1$
We can now apply Gronwall inequality since
\paragraph{Verlet}
For Verlet, using Lemmas 4, 5, 7, 9, we get, 
\begin{eqnarray*}
e^{n+1} & \leq & C(h^2+\Delta t^3 + (\Delta t^2+\Delta t)(h^2+ h^2\Delta t \Delta t^2 + e^n)+e^n,\nonumber\\
	     & \leq & (1+C\Delta t)e^n+ C(h^2+\Delta t^3 + (\Delta t^2+\Delta t)(h^2+\Delta t^2)),
\end{eqnarray*}
So that 
$$
e^n \leq \exp(C'T)e^0 + C(h^2+\Delta t^2 + \frac{h^2}{\Delta t}),
$$
which is what was expected. Taking $\Delta t = Ch^{\alpha}$, we find the best global order with $\alpha = \frac{2}{3}$ being $\frac{4}{3}$.

\paragraph{CK2}
For $CK2$, using Lemmas 4, 5, 8, 10, we get 
\begin{eqnarray*}
e^{n+1} & \leq & C(h^2+\Delta t^3+ (\frac{\Delta t^2}{\Delta x}+\Delta t)(h^2+ e^n))+e^n,\nonumber\\
	& \leq &  (1+C(\Delta t+\frac{\Delta t^2}{\Delta x})e^n+ C(h^2+\Delta t^3 + \Delta t \,h^2+ \frac{\Delta t^2}{\Delta x}h^2),
\end{eqnarray*}
Here, assuming $\Delta t \leq \Delta x$, we have :
$$
e^n \leq \exp(C'T)e^0 + C(h^2+\Delta t^2 + \frac{h^2}{\Delta t}),
$$
If you want to look for the best global order here, you find the same result as in Verlet, nevertheless, this cannot fit with the other assumption $\Delta t \leq \Delta x$. Therefore, the only way is to take $\Delta t = \Delta x$, and the global order is 1.
\paragraph{CK3}
For $CK3$, using the same lemmas as for $CK2$:
\begin{eqnarray*}
e^{n+1} & \leq & C(h^2+\Delta t^4 + e^n(\Delta t + \frac{\Delta t^2}{\Delta x} + \frac{\Delta t^3}{\Delta x^2})+\Delta t h^2 + \Delta t^2 \frac{h^2}{\Delta x}+\Delta t^3 \frac{h^2}{\Delta x^2})+e^n,\nonumber\\
	& \leq & e^n(1+C(\Delta t + \frac{\Delta t^2}{\Delta x} + \frac{\Delta t^3}{\Delta x^2}))+C(h^2 + \Delta t h^2+ \frac{\Delta t^2}{\Delta x}h^2+\frac{\Delta t^3}{\Delta x^2}h^2 + \Delta t^4).
\end{eqnarray*}
Assuming again $\Delta t \leq \Delta x$, we have 
$$
e^n \leq \exp(C'T)e^0 + C(\Delta t^3+h^2+\frac{h^2}{\Delta t}).
$$
The same remark as with $CK2$ is still valid. We can see that we are limited because of the terms $\frac{h^2}{\Delta t}$ and $\frac{\Delta t}{\Delta x}$. In order to be able to reach higher orders, we would have to use splines of superior degrees $m>1$ to get terms in $\frac{h^{m+1}}{\Delta t}$ like in the other proofs of convergence, for example: (\cite{mehrenberger,bostan})



 
\section{Numerical results}

In order to validate our new schemes, we have tested them on two standard test cases of plasma physics, the two stream instability and the bump on tail instability. We also compared them to the classical and knowledgeably robust Verlet scheme. Notice that because of the diffusitivity of linear splines, we have used cubic splines for the distribution function.

\begin{figure}
\begin{tabular}{cc}
\includegraphics[width=8cm]{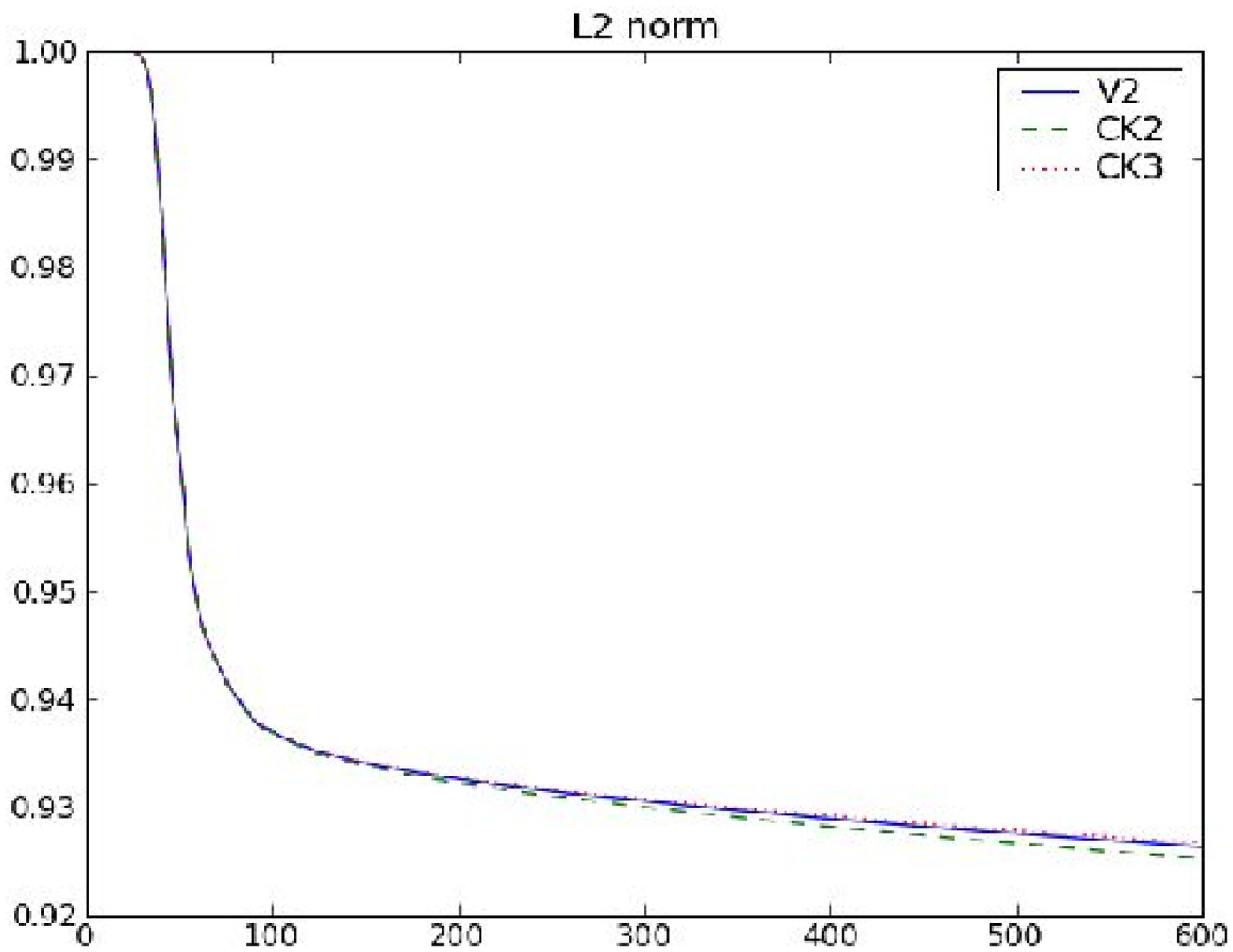} & \includegraphics[width=8cm]{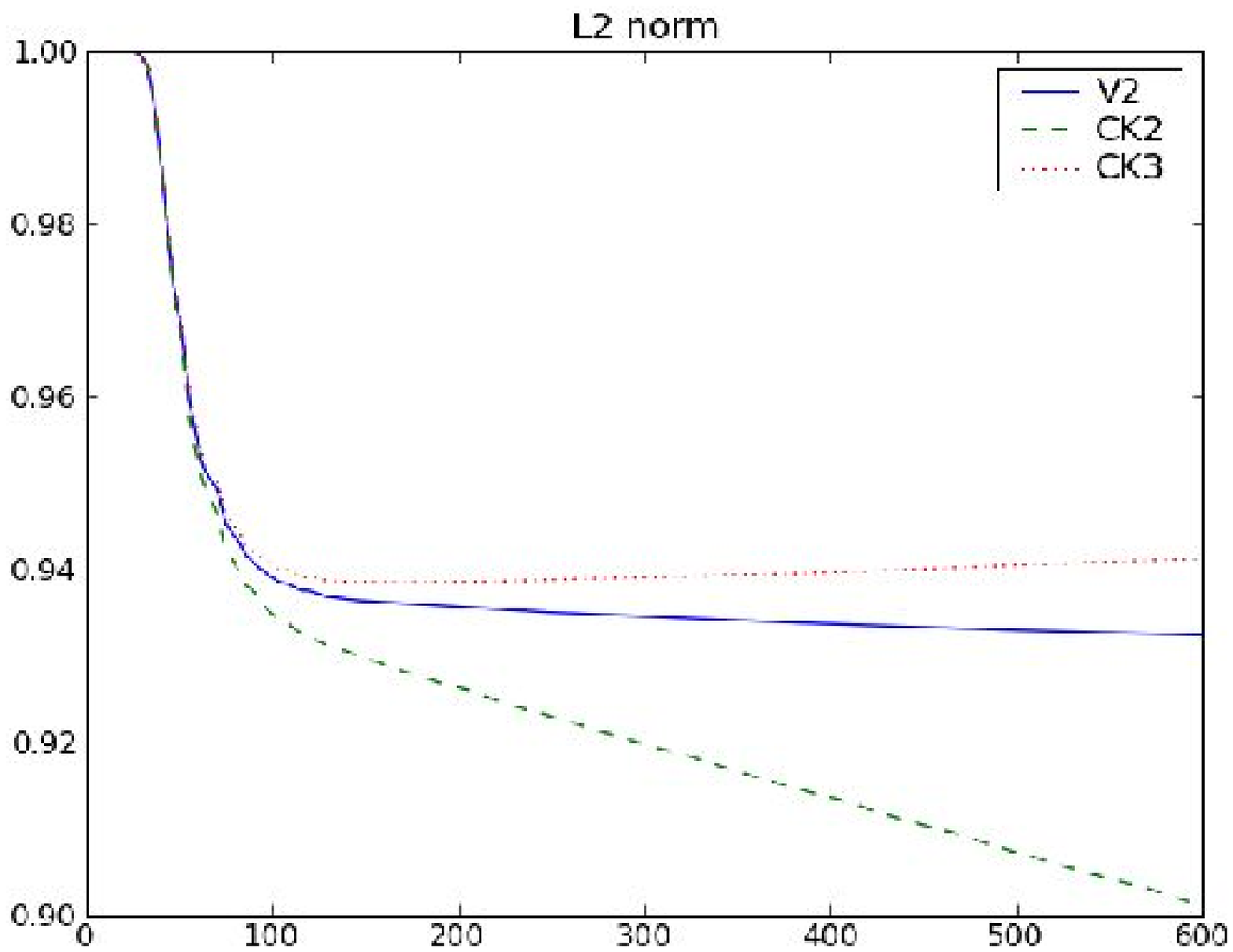}\\
\includegraphics[width=8cm]{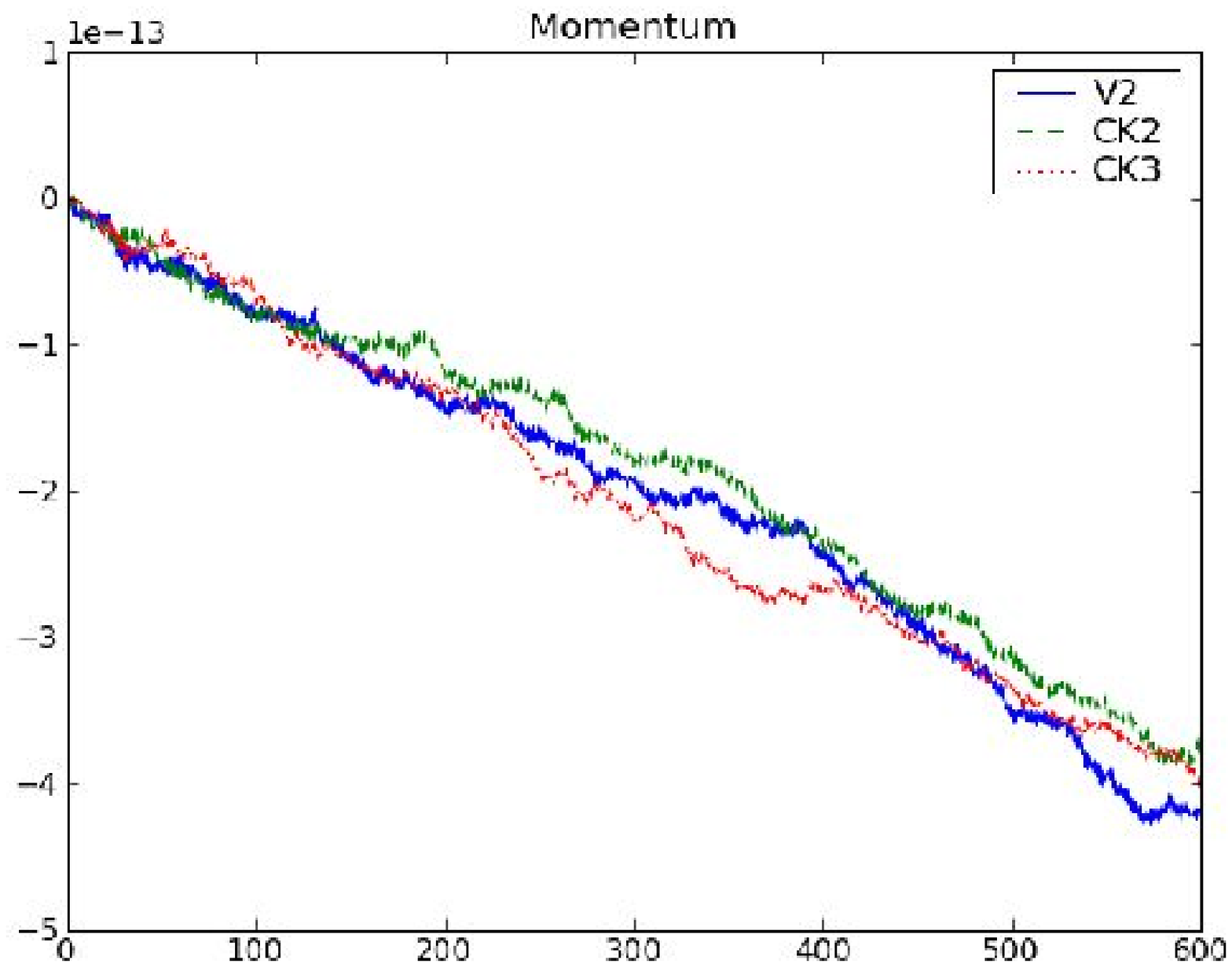} &\includegraphics[width=8cm]{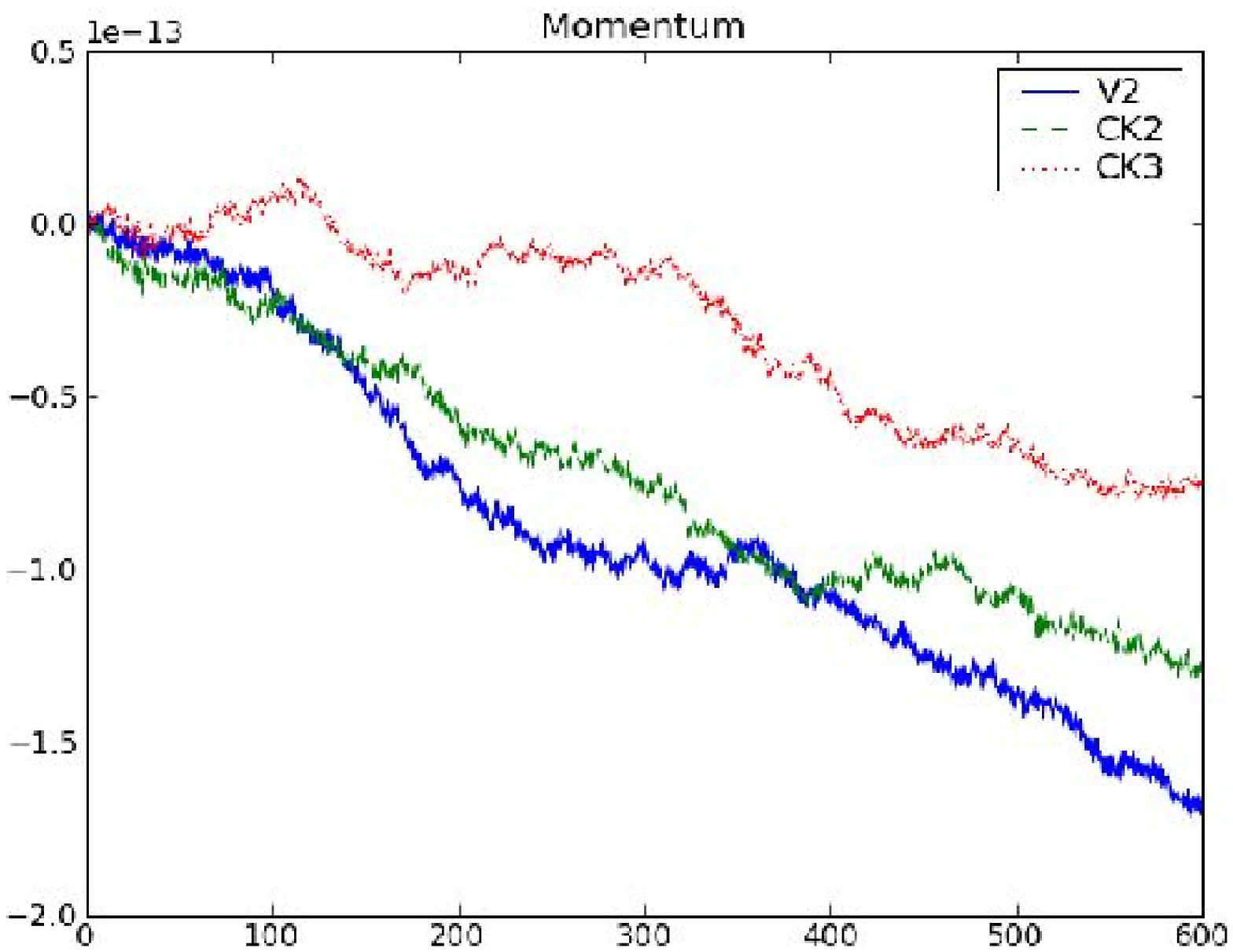}\\
\includegraphics[width=8cm]{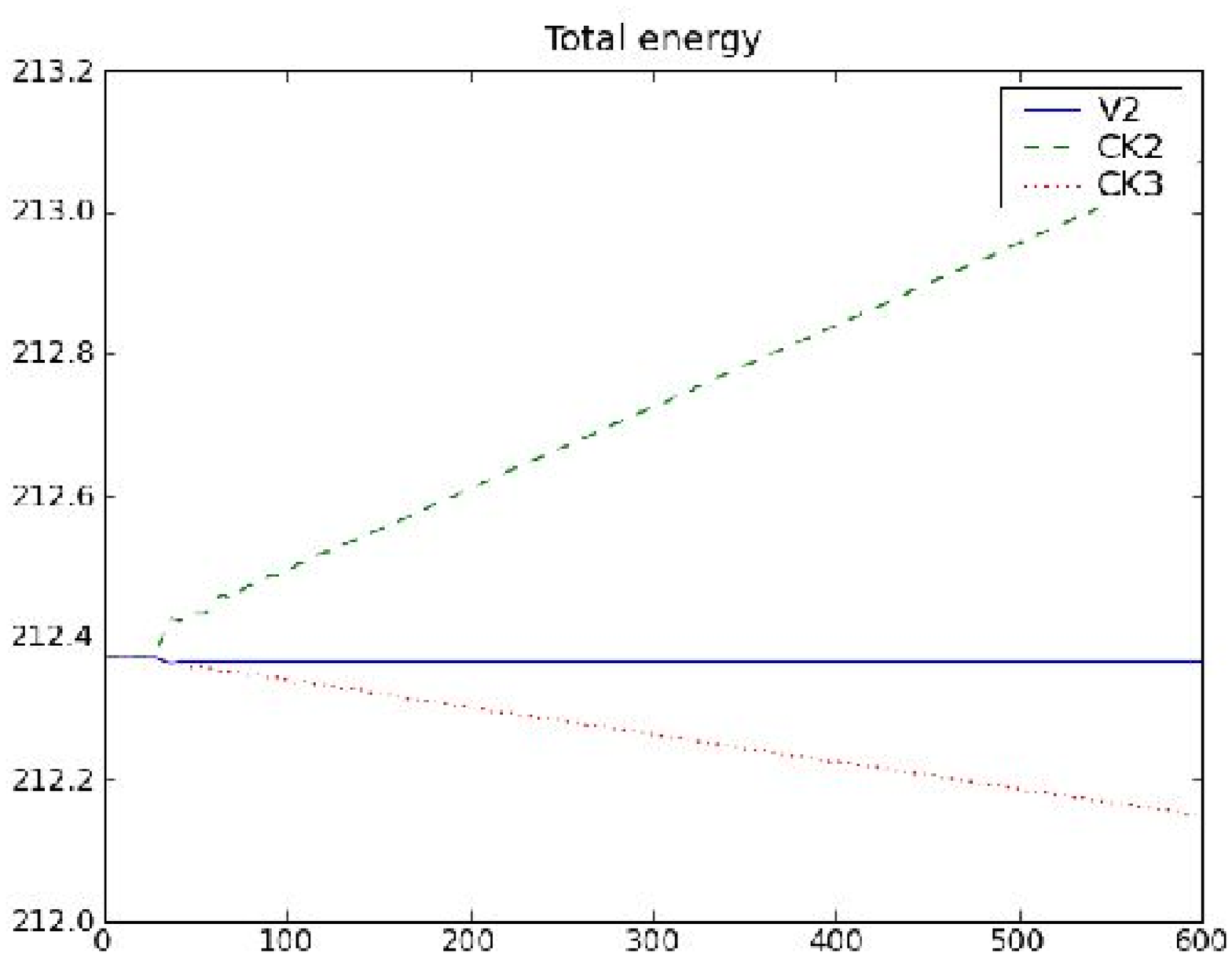} &\includegraphics[width=8cm]{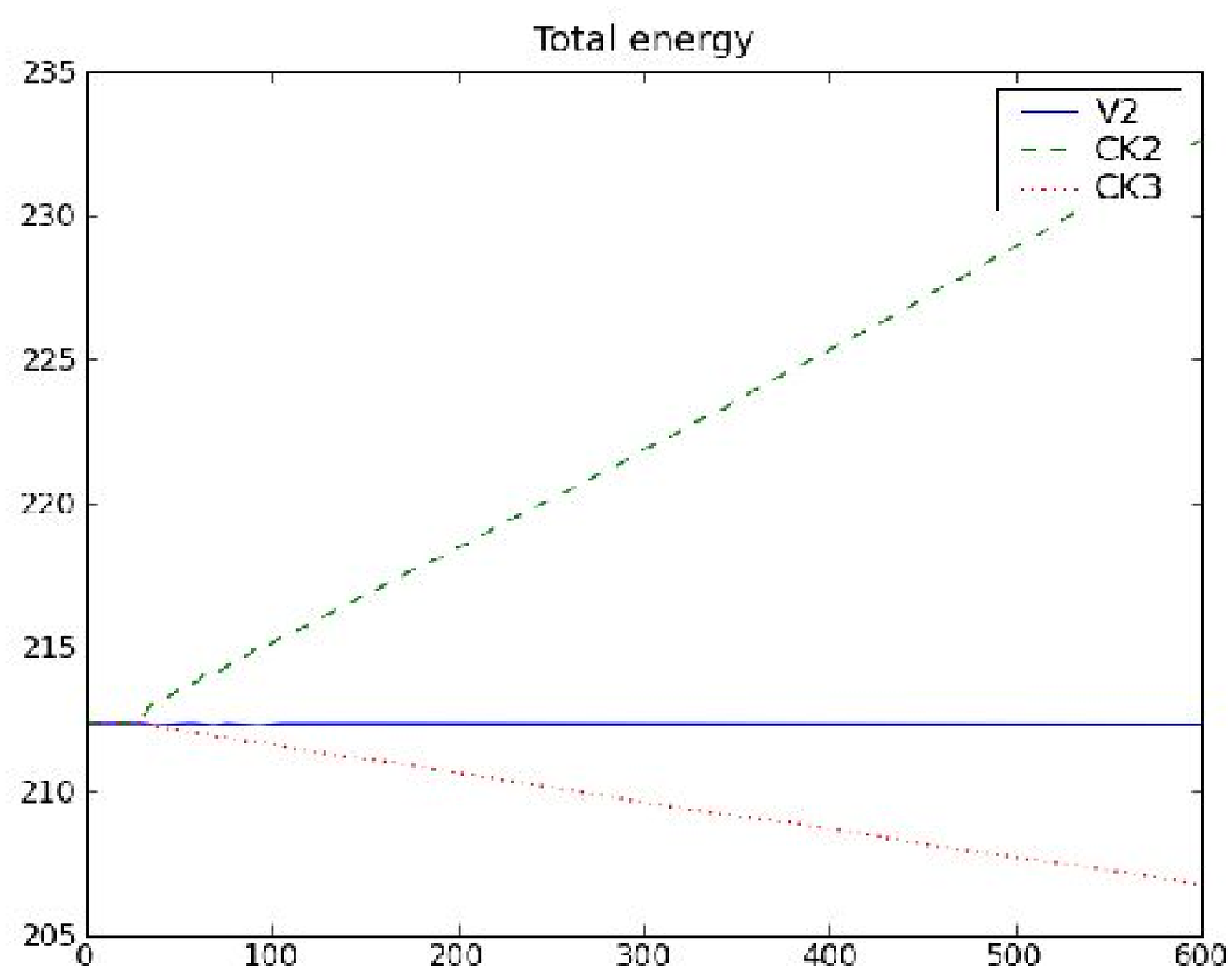}  
\end{tabular}
\caption{\label{fig:tsi}Two stream instability:}
\end{figure}

For the two stream instability, the initial condition is given by
$$
f_0(x, v) =\frac{1}{\sqrt{2\pi}}e^{-v^2/2}v^2[1-\alpha\cos(k x)], 
$$ 
with $k=0.2$ and $\alpha=0.05$. The computational domain is 
$[0,2\pi/k]\times [-9,9]$ which is sampled by $N_x=N_v=128$ points. 
We used a time step $\Delta t=0.1$ in the results of the left hand side of Figure \ref{fig:tsi} and of 
$\Delta t=0.3$ on the right-hand side of the figure.
We display the 
$L^2$ norm which reveals the dissipation of the scheme, the total momentum and the total energy.
All of those are conserved in the continuous Vlasov-Poisson system. We do not display the number of particles which is conserved with an even better accuracy than the momentum. The momentum is exactly conserved by the scheme and up to about $10^{-13}$ in the simulation. This is due to roundoff errors and the truncation of the velocity space. The $L^2$ norm cannot be exactly conserved by any scheme using a phase space grid as soon as the grid does not resolve anymore the filaments.
The Verlet scheme is our reference scheme here, and we observe that the results obtained with the CK schemes are very close, especially for the smallest time step. Moreover conservation properties are better for the third order CK3 than for the second order CK2.

\begin{figure}
\begin{tabular}{cc}
\includegraphics[width=8cm]{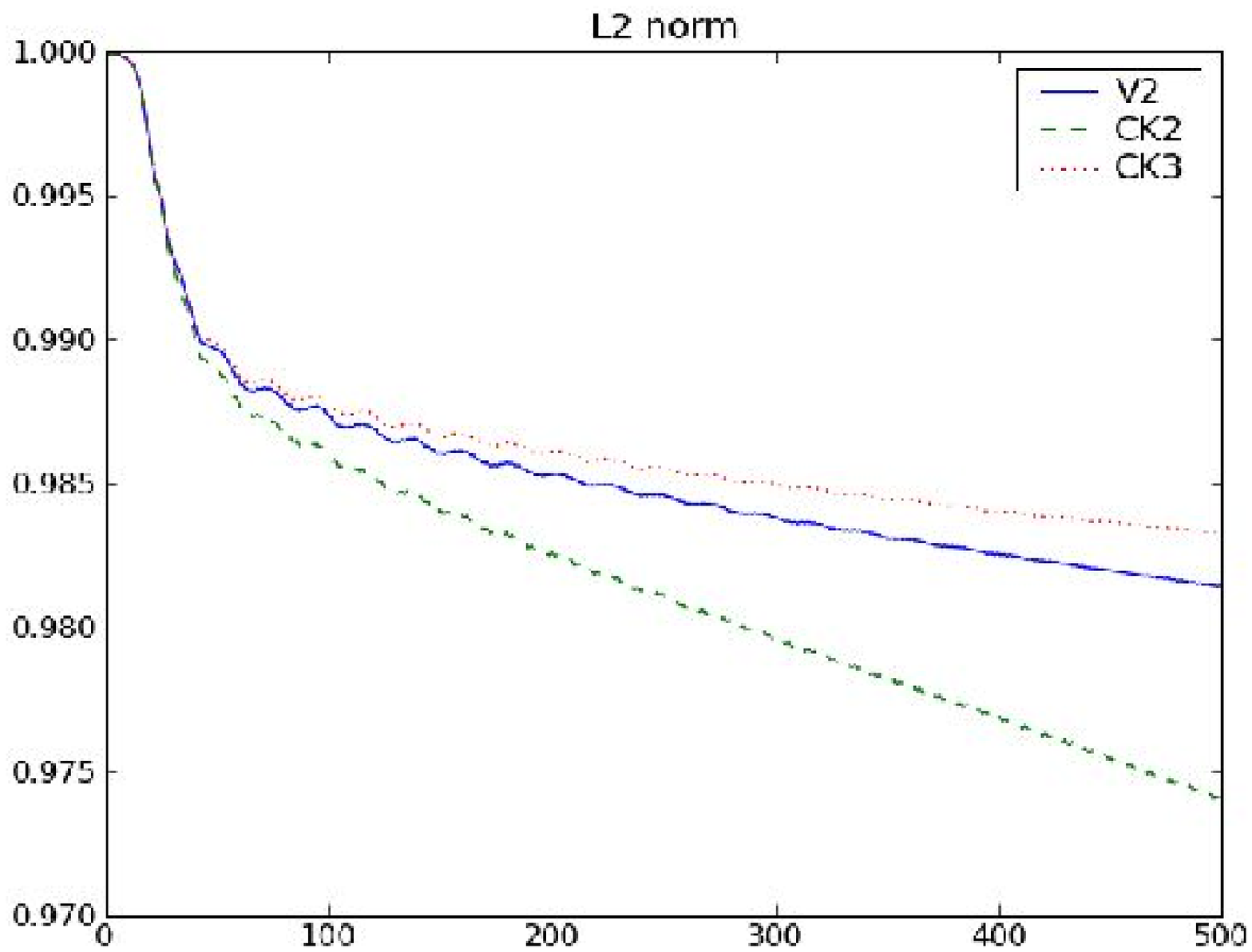} & \includegraphics[width=8cm]{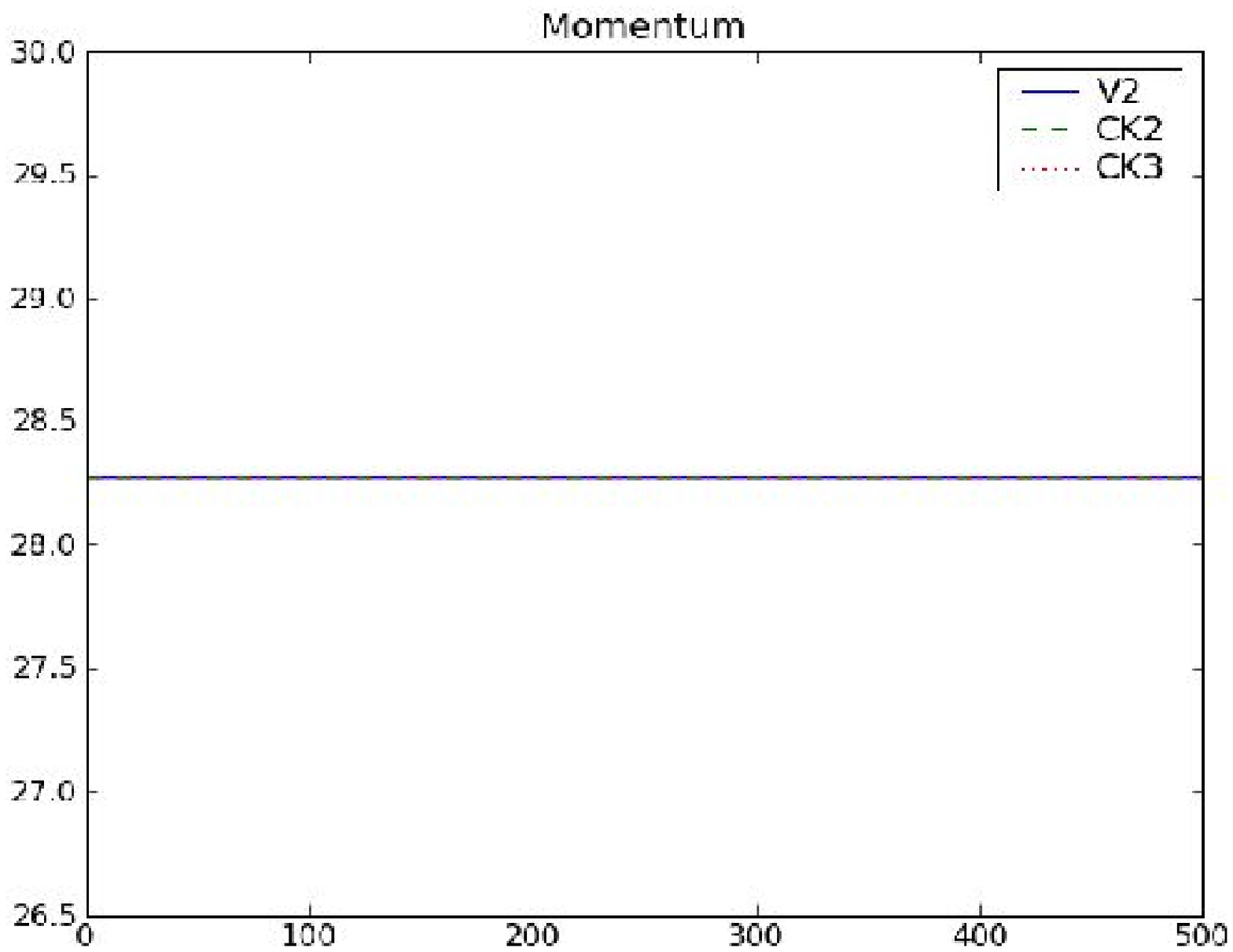}\\
\includegraphics[width=8cm]{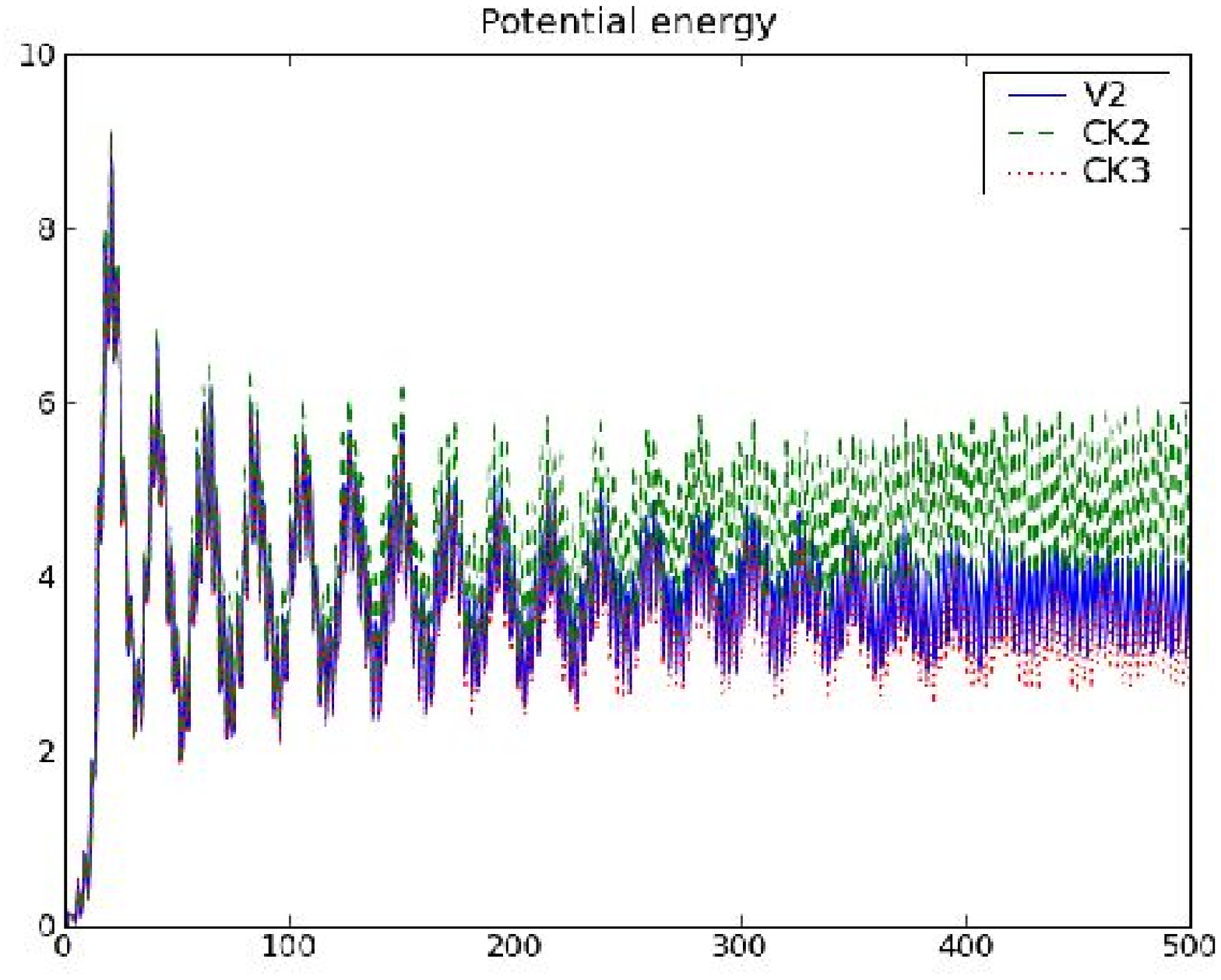} &\includegraphics[width=8cm]{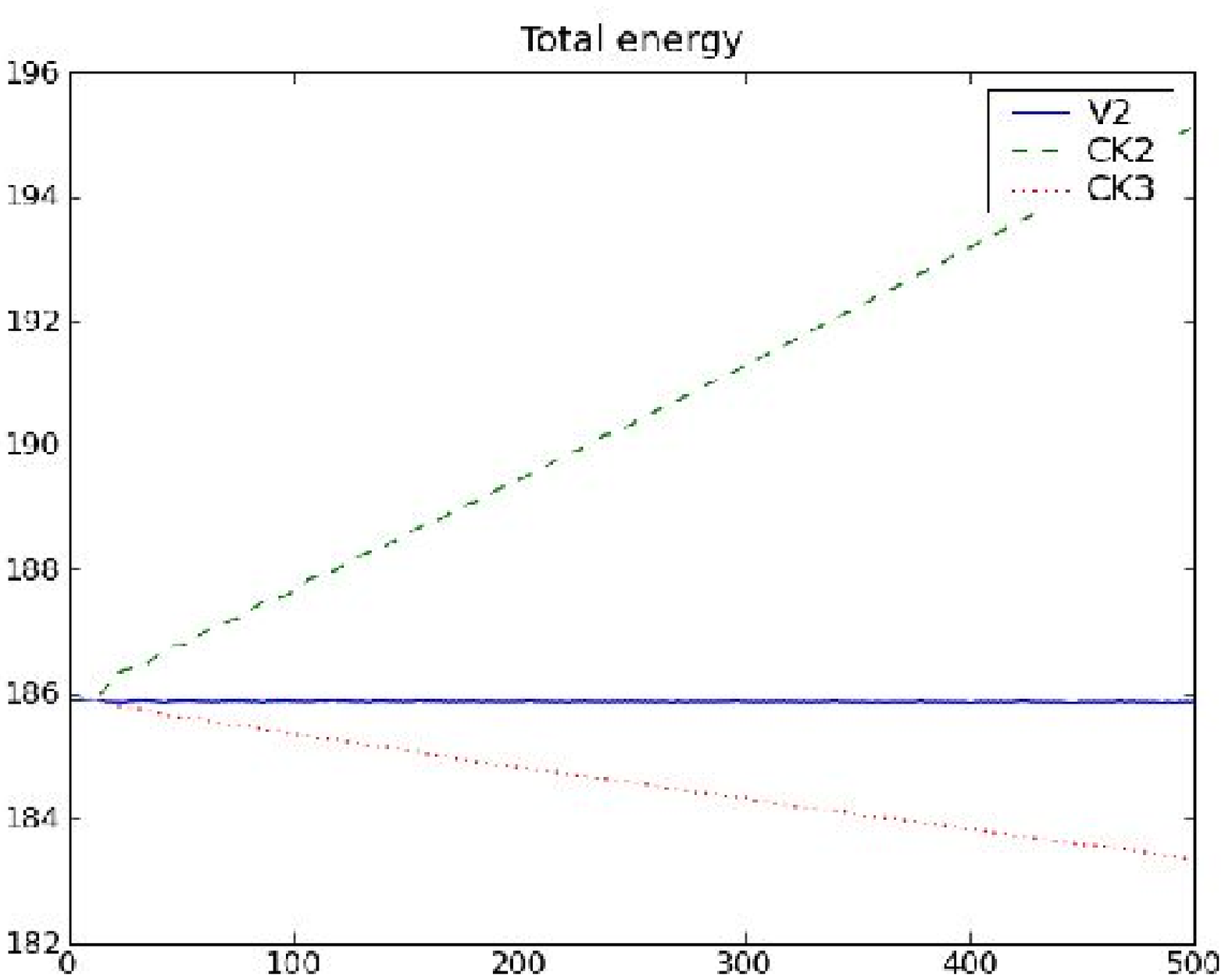}  
\end{tabular}
\caption{\label{fig:bot}Bump on tail}
\end{figure}

For the bump-on-tail instability test case 
 the initial condition writes 
$$
f_0(x, v) = \tilde{f}(v)[1+\alpha\cos(k x)], 
$$ 
with 
$$
\tilde{f}(v)=n_p \exp(-v^2/2) + n_b \exp\left(-\frac{|v-u|^2}{2v^2_t}\right)
$$
on the interval $[0, 20\pi]$, with periodic conditions in space. 
The initial condition $f_0$ is a Maxwellian distribution function which has a bump on 
the Maxwell distribution tail; the parameters of this bump are the following 
$$
n_p=\frac{9}{10(2\pi)^{1/2}}, n_b=\frac{2}{10(2\pi)^{1/2}}, u=4.5, v_t=0.5, 
$$
whereas the numerical parameters are $N_x=128, N_v=128, v_{max}=9, \Delta t=0.2$. 
The results are displayed in Figure \ref{fig:bot}. Here the momentum does not vanish, so that the results are not polluted by roundoff errors and the momentum is exactly conserved.
The conclusion for the $L^2$ norm and the total energy is the same as in the Two Stream Instability test case. The potential or electric energy is a classical diagnostic for the bump on tail instability. The oscillations go on for a long time with all three time schemes, even though there is a slight energy increase for the CK2 scheme.

\section{Conclusion}

In this paper, the proof of a $L^1$ convergence has been reached for linear spline interpolation. The originality, except from the choice of the $L^1$ norm is that the convergence has been reached for a non split method. In this paper, the computation of the characteristics has been made with the Verlet algorithm, or with a CK procedure, but the proof can be adapted to other algorithms such as Runge Kutta of any order. There remains for the moment some problems using splines of higher orders, especially concerning stability. This prevents us from reaching real high order algorithms.  Numerical experiments that can be seen in \cite{respaud}, and confirmed here seem to prove that the method is also stable and convergent for cubic splines. Nevertheless, there remains a problem to preserve the $l^1$ norm of the coefficients $\omega_{i,j}$, since some of them can become non positive in the solving of the linear system with splines of degree higher than 2. Another way of tackling the problem will probably be needed.

\cleardoublepage


\begin{thebibliography}{1}


\bibitem{mehrenberger}
{\sc N. Besse, M. Mehrenberger}, 
\textit{Convergence of classes of high order semi-Lagrangian schemes for the Vlasov-Poisson system}, 
Math. Comput., {\bf 77}, pp. 93--123, (2008). 

\bibitem{bostan}
{\sc M. Bostan, N. Crouseilles},
\textit{Convergence of a semi-Lagrangian scheme for the reduced Vlasov-Maxwell system for laser-plasma interaction},
Numer. Math. {\bf 112}, pp. 169--195, (2009).

\bibitem{bouchut}
{\sc F. Bouchut, F. Golse, M. Pulvirenti},
\textit{Kinetic equations and asymptotic theory,}
Series in applied Math. P.G Ciarlet and P.L Lions (Eds), Gauthier Villars (2008).


\bibitem{birdsall}
{\sc C.K. Birdsall, A.B. Langdon}, 
\textit{Plasma Physics via Computer Simulation}, 
Inst. of Phys. Publishing, Bristol/Philadelphia, 1991. 

\bibitem{vecil}
{\sc J.-A. Carillo, F. Vecil}, 
\textit{Non oscillatory interpolation methods applied to Vlasov-based models}, SIAM J. Sci. Comput. {\bf 29}, pp. 1179--1206, (2007). 

\bibitem{cheng}
{\sc C. Z. Cheng, G. Knorr}, 
\textit{The integration of the Vlasov equation in configuration space}, 
J. Comput. Phys, {\bf 22}, pp. 330--3351, (1976). 

\bibitem{respaud}
{\sc N. Crouseilles, T. Respaud, E. Sonnendr\"ucker},
\textit{A forward semi-Lagrangian method for the numerical solution of the Vlasov equation},
Comput. Phys. Comm., {\bf 180} (10), pp. 1730--1745, (2009).

\bibitem{Cotter}
{\sc C.J. Cotter, J. Frank, S. Reich}
\textit{The remapped particle-mesh semi-Lagrangian advection scheme,}
Q. J. Meteorol. Soc., {\bf 133}, pp. 251--260, (2007).

\bibitem{Cottet}
{\sc G.-H Cottet, P.-A Raviart},
\textit{Particle methods for the one-dimensional Vlasov-Poisson equations,}
Siam J. Numer. anal. {\bf 21}, pp. 52--75, (1984).

\bibitem{Despres1}
{\sc B. Despr\'es},
\textit{Finite volume transport Schemes,}
Numerische Mathematik {\bf 108}, pp.529--556, (2008).

\bibitem{filbet1}
{\sc F. Filbet, E. Sonnendr\"{u}cker, P. Bertrand}, 
\textit{Conservative numerical schemes for the Vlasov equation},
J. Comput. Phys.,  {\bf 172}, pp. 166--187, (2001).

\bibitem{filbet2}
{\sc F. Filbet, E. Sonnendr\"{u}cker}, 
\textit{Comparison of  Eulerian Vlasov solvers},
Comput. Phys. Comm.,  {\bf 151}, pp. 247--266, (2003).

\bibitem{Glassey}
{\sc R.T Glassey},
\textit{The Cauchy problem in kinetic theory},
SIAM, Philadelphia (1996).


\bibitem{virginie}
{\sc V. Grandgirard, M. Brunetti, P. Bertrand, N. Besse, X. Garbet, P. Ghendrih, G. Manfredi, Y. Sarrazin, O. Sauter, E. Sonnendr\"{u}cker, J. Vaclavik, L. Villard}, 
\textit{A drift-kinetic semi-Lagrangian 4D code for ion turbulence simulation}, 
J. Comput. Phys., {\bf 217}, pp. 395--423, (2006). 

\bibitem{Reich} 
{\sc S. Reich}, 
\textit{An explicit and conservative remapping strategy for semi-Lagrangian advection}, 
Atmospheric Science Letters {\bf 8}, pp. 58--63, (2007). 

\bibitem{staniforth}
{\sc A. Staniforth, J. Cot\'e}, 
\textit{Semi-Lagrangian integration schemes for atmospheric models - A review}, 
Mon. Weather Rev. {\bf 119}, pp. 2206--2223, (1991). 

\bibitem{sonnen}
{\sc E. Sonnendr\"{u}cker, J. Roche, P. Bertrand, A. Ghizzo}
\textit{The semi-Lagrangian method for the numerical resolution of the Vlasov equation,}
J. Comput. Phys.,  {\bf 149}, pp. 201--220, (1999).

\bibitem{ZerWoS05}
{\sc M. Zerroukat, N. Wood, A. Staniforth},
\textit{ A monotonic and positive-definite filter for a Semi-Lagrangian Inherently Conserving and Efficient (SLICE) scheme}, Q.J.R. Meteorol. Soc., {\bf 131}, pp. 2923-2936, (2005).

\bibitem{ZerWoS06}
{\sc M. Zerroukat, N. Wood, A. Staniforth},
\textit{The Parabolic Spline Method (PSM) for conservative transport problems},
Int. J. Numer. Meth. Fluids, {\bf  51}, pp. 1297--1318, (2006).

\end{thebibliography}
\end{document}